\documentstyle{amsppt}
\pageheight{19cm}
\magnification=1200
\baselineskip=18pt
\nologo
\TagsOnRight

\document

\define \ts{\thinspace}

\define \wh{\widehat}
\define \ot{\otimes}
\define \sdet{\text{{\rm sdet}\ts}}
\define \tr{\text{{\rm tr}\ts}}
\define \oa{{\frak o}}
\define \Sym{{\frak S}}
\define \spa{{\frak {sp}}}
\define \U{{\operatorname {U}}}
\define \A{{\operatorname {A}}}
\define \G{{\operatorname {G}}}
\define \GL{{\operatorname {GL}}}
\def\SO{{\operatorname {SO}}}
\def\SP{{\operatorname {Sp}}}
\def\P{{\operatorname {P}}}
\define \Q{{\operatorname {\Lambda^*}}}
\define \M{{\operatorname {M^*}}}
\define \Qc{{\operatorname {\Lambda_c^*}}}
\define \Mc{{\operatorname {M_c^*}}}
\define \Z{{\operatorname {Z}}}
\define \I{{\operatorname {I}}}
\define \J{{\operatorname {J}}}
\define \LL{{\operatorname {L}}}
\define \Y{{\operatorname {Y}}}
\define \C{{\Bbb C}}
\define \ZZ{{\Bbb Z}}
\define \gr{\text{\rm gr}\ts}
\define \sgn{\text{\rm sgn}\ts}
\define \qdet{\text{{\rm qdet}\ts}}
\define \End{\text{{\rm End}\ts}}
\define \T{{\Cal T}}
\define \hra{\hookrightarrow}
\define \Proof{\noindent {\it Proof.}\ \ \ }
\define \g{{\frak g}}
\define \h{{\frak h}}
\define \n{{\frak n}}
\define \gti{\widetilde{\frak{g}}}
\define \hti{\widetilde{\frak{h}}}
\define \nti{\widetilde{\frak{n}}}
\def\l{{\lambda}}
\define \gl{{\frak{gl}}}
\define \gln{{\frak{gl}(n)}}

\heading{\bf CENTRALIZER CONSTRUCTION FOR TWISTED YANGIANS} \endheading

\bigskip
\bigskip
\heading{Alexander Molev and Grigori Olshanski}
\endheading 
\bigskip
\bigskip
\heading{\bf December 1997}\endheading

\bigskip
\bigskip
\noindent
{\bf Abstract}\newline
For each of the classical Lie algebras
$\g(n)=\oa(2n+1),\spa(2n),\oa(2n)$ of type $B$, $C$, $D$ we consider
the centralizer of the subalgebra $\g(n-m)$ in the universal
enveloping algebra $\U(\g(n))$. We show that
the $n$th centralizer algebra can be naturally
projected onto the $(n-1)$th one, so that one can form the projective
limit of the centralizer algebras as $n\to\infty$ with $m$ fixed.
The main result of the paper is a precise description of this limit
(or stable) centralizer algebra, denoted by $\A_m$. We explicitly
construct an algebra isomorphism $\A_m=\Z\ot \Y_m$, where $\Z$ is a
commutative algebra and $\Y_m$ is the so-called twisted Yangian
associated to the rank $m$ classical Lie algebra of type $B$, $C$, or
$D$. The algebra $\Z$ may be viewed as the algebra of virtual Laplace
operators; it is isomorphic to the algebra of polynomials with
countably many indeterminates. The twisted Yangian $\Y_m$ (and hence
the algebra $\A_m$) can be described in terms of a system of generators
with quadratic and linear defining relations which are conveniently
presented in $R$-matrix form involving the so-called reflection
equation. This extends the earlier work on the type $A$ case
by the second author. 


\bigskip
\bigskip
\noindent
{\bf Postal addresses:}\newline
Centre for Mathematics and its Applications\newline
Australian National University\newline
Canberra, ACT 0200, Australia\newline
(E-mail:\ molev\@pell.anu.edu.au)
\medskip
\noindent
Dobrushin Mathematics Laboratory\newline
Institute for Problems of Information Transmission\newline
Bolshoy Karetny 19, Moscow 101447, GSP-4, Russia\newline
(E-mail:\ olsh\@ippi.ac.msk.su,\  bv\@glasnet.ru)

\newpage

\noindent
{\bf 0. Introduction}
\bigskip

Let $\g$ be a complex reductive Lie algebra, let $\g'\subset\g$ be a
reductive subalgebra, and let $\Z(\g,\g')$ denote the centralizer of $\g'$
in the universal enveloping algebra $\U(\g)$.
The centralizer algebra $\Z(\g,\g')$ appears, for example, in the following
situation. Assume $V$ is an irreducible $\g$-module which decomposes
(under the action of $\g'$) into a direct sum of irreducible
finite-dimensional $\g'$-modules $W_\lambda$ with certain
multiplicities $\text{mult}_\lambda$ (here $\{\lambda\}$ is the
set of dominant highest weights for $\g'$ that occur in $V$);
then this decomposition can be written as
$$
V\simeq \sum_\lambda U_\lambda\otimes W_\lambda,
\qquad \text{dim}\ts U_\lambda=\text{mult}_\lambda, 
$$
where, for each $\lambda$,
$$
U_\lambda=\text{Hom}_{\g'}(W_\lambda,V)
$$
is an irreducible $\Z(\g,\g')$-module; see, e.g., Dixmier [D], Section 9.1.

For some special couples $(\g,\g')$ the centralizer algebra turns out
to be commutative and so $\text{mult}_\lambda\equiv1$. This holds,
for example,
for the couples $(\frak{gl}(N), \frak{gl}(N-1))$
and $(\frak{o}(N),\frak{o}(N-1))$, but even for the allied couple
$(\frak{sp}(2n),\frak{sp}(2n-2))$ the centralizer algebra is
noncommutative and its structure seems to be complicated.
It is believable that an understanding of the algebra
$\Z(\frak{sp}(2n),\frak{sp}(2n-2))$ and its representations could lead
to a solution of an old problem -- the construction of an orthonormal
Gelfand--Tsetlin-type basis for representations of $\frak{sp}(2n)$.

In the present paper we investigate certain series of couples $(\g,\g')$
which are indexed by two parameters $m<n$:
$$
\alignedat 4
&\text{type $A$} & \qquad &\text{type $B$} & \qquad &\text{type $C$} &
\qquad &\text{type $D$}\\
&\frak{gl}(n) & &\frak{o}(2n+1) & &\frak{sp}(2n) & &\frak{o}(2n)\\
&\frak{gl}(n-m) & &\frak{o}(2(n-m)+1) & &\frak{sp}(2(n-m)) & 
&\frak{o}(2(n-m)),
\endalignedat \tag 0.1
$$
and we study the `stable structure' of the centralizer $\Z(\g,\g')$ as
$n\to\infty$ with $m$ fixed.
In more detail, write $\g=\g(n)$, $\g'=\g_m(n)$ and abbreviate
$$
\A_m(n)=\Z(\g(n),\g_m(n)).
$$
It turns out that for any fixed $m$ there exist natural projections
$$
\pi_n: \A_m(n)\to \A_m(n-1)
$$ 
which are algebra morphisms and which preserve filtration (induced by the
standard filtration of enveloping algebras), so that one can form the
projective limit
$$
\A_m=\underset{\longleftarrow}\to{\lim}\ts (\A_m(n),\pi_n), 
\qquad n\to\infty
$$
in the category of filtered algebras. 
We call this the {\it centralizer construction\/}.

The main result of this paper is that each of the limit algebras $\A_m$ (in
contrast to the centralizers $\A_m(n)$) admits a very precise
description. Namely, we find for $\A_m$ a system of generators with
{\it quadratic\/} and {\it linear\/} 
defining relations which are conveniently written using the $R$-{\it
matrix formalism\/}. Moreover, this kind of presentation of the
algebra $\A_m$ makes it possible to study its representations (and hence
representations of the centralizer algebras).

For couples $(\g,\g')$ of type $A$ this was done earlier in
Olshanski [O1], [O2]. In that case 
$$
\A_m=\A_0\otimes \Y(\frak{gl}(m)),
$$  
where $\A_0$ is a commutative algebra (isomorphic to the ring of
polynomials with countably many indeterminates) and $\Y(\frak{gl}(m))$
is the {\it Yangian\/} for the Lie algebra $\frak{gl}(m)$.
The algebra $\Y(\frak{gl}(m))$ first appeared in the works of
L.~D.~Faddeev's school on the Yang--Baxter equation; see, e.g.,
[TF]. One starts with the `ternary relation'
$$
R(u-v)(T(u)\otimes 1)(1\otimes T(v))=
(1\otimes T(v)) (T(u)\otimes 1) R(u-v),\tag 0.2
$$
where $u$, $v$ are formal parameters, $R(u)$ is {\it Yang's $R$-matrix\/},
$$
R(u)=1\otimes1-u^{-1} \sum_{i,j=1}^m E_{ij}\otimes E_{ji}\in 
\text{End}(\Bbb C^m\otimes \Bbb C^m)(u),
$$
the $E_{ij}$ are matrix units, and $T(u)$ is a matrix--valued formal
series in $u^{-1}$:
$$
T(u)=\left(t_{ij}(u)\right), 
\qquad t_{ij}(u)=\delta_{ij}+\sum_{k=1}^\infty t_{ij}^{(k)}u^{-k},
\quad 1\le i,j\le m.
$$
Then the relation \thetag{0.2} implies a system of quadratic
relations on the symbols $t_{ij}^{(k)}$ which is just the system of
defining relations for $\Y(\frak{gl}(m))$.
Note that $\Y(\frak{gl}(m))$ is a Hopf algebra and that it can be
viewed as a deformation of the enveloping algebra $\U(\frak{gl}(m)[x])$, 
where
$$
\frak{gl}(m)[x]=\frak{gl}(m)+\frak{gl}(m)x+\frak{gl}(m)x^2+\dots
$$
is the Lie algebra of $\frak{gl}(m)$-valued polynomials (a 
polynomial current Lie algebra). 

V.~G.~Drinfeld defined Yangians
corresponding to arbitrary simple Lie algebras [D1], [D2]. The Yangians
form an important family of quantum groups and have a rich representation
theory; see [CP1], [CP2], [C2], [D2].

One could expect the algebras $\A_m$ associated with type 
$B, C, D$ couples $(\g,\g')$ in \thetag{0.1} to be connected
with orthogonal 
or symplectic Drinfeld's Yangians. However, this is not at all true,
and we have the following result (announced in [O3]).

\proclaim{\bf Main Theorem}
For the couples $(\g,\g')$ of type $B, C, D$ one has
$$
\A_m=\Z\otimes \Y^{\pm}(M),
$$
where $\Z$ is a commutative algebra, $M=2m$ or $M=2m+1$
{\rm(}for types $C,D$ or $B$, respectively{\rm)}, and
$\Y^{\pm}(M)$ is 
a certain `Yangian-type' quadratic algebra which can be realized as 
a one-sided Hopf ideal in the Yangian $\Y(\frak{gl}(M))$.
\endproclaim

It should be emphasized that the algebras $\Y^{\pm}(M)$ are not Hopf
algebras. We call these new objects {\it twisted Yangians\/}
(the sign `$+$' is taken in the orthogonal case, the sign `$-$'
is taken in the symplectic case). 

Like the Yangian $\Y(\frak{gl}(m))$, the algebra $\Y^{\pm}(M)$ is a
deformation of an enveloping algebra; the 
corresponding Lie algebra is an involutive subalgebra of 
$\frak{gl}(M)[x]$: 
$$
\{f\in\frak{gl}(M)[x]\mid f(-x)=-(f(x))^t\}, 
$$
where $t=t_{\pm}$ stands for the matrix transposition
which corresponds to a symmetric or alternating form in $\Bbb C^M$.
This is a `twisted' polynomial current Lie algebra, which explains
the term `twisted Yangian'. 

There are two kinds of defining relations for generators of
$\Y^{\pm}(M)$: quadratic and linear, and both are conveniently
written in an $R$-matrix form.
The presentation of the quadratic relations is similar to \thetag{0.2} but
more complicated:
$$
R(u-v)(S(u)\otimes 1)R^t(-u-v)(1\otimes S(v))=
(1\otimes S(v))R^t(-u-v)(S(u)\otimes 1) R(u-v),
\tag 0.3
$$
where 
$$
R^t(u)=(\text{id}\otimes t)(R(u))
=(t\otimes\text{id})(R(u)).
$$

As with the algebra $\Y(\frak{gl}(m))$,
such a modification of the relation \thetag{0.2} also first appeared
in mathematical physics (see Cherednik [C1] and Sklyanin [S]); 
nowadays it is called the {\it reflection
equation\/}, see, e.g., [KK], [KS], [KJC].  

The system of linear relations for the generators of $\Y^{\pm}(M)$ 
is written as the following {\it symmetry relation\/} 
on the matrix $S$:
$$
S^t(-u)=S(u)\pm\frac{S(u)-S(-u)}{2u}.\tag 0.4
$$ 
The explanation of this relation is rather simple: it reflects the
fact that the matrices $X\in\g(m)$ are antisymmetric with respect to the
transposition $t$.

Let us describe one more aspect of the centralizer construction.
Existence of natural embeddings $\g(n)\to\g(n+1)$ allows us to form
the inductive limit Lie algebra 
$$
\g(\infty)=\underset{\longrightarrow}\to{\lim}\ts\g(n),\qquad n\to\infty
$$
which is one of the algebras
$$
\alignat 4
&\text{type $A$} & \qquad &\text{type $B$} & \qquad &\text{type $C$} &
\qquad &\text{type $D$}\\
&\frak{gl}(\infty) & &\frak{o}(2\infty+1) & 
&\frak{sp}(2\infty) & &\frak{o}(2\infty).
\endalignat
$$
(Note that 
an appropriate modification of the centralizer construction
can be also applied to another version of the algebra
$\frak{gl}(\infty)$, denoted by {$\frak{gl}(2\infty+1)$}, 
see Remark 2.23.)

By its very definition, $\A_m(n)$ is a subalgebra
of $\A_{m+1}(n)$, which leads to an algebra embedding $\A_m\to \A_{m+1}$,
so that we can form the inductive limit algebra
$$
\A=\underset{\longrightarrow}\to{\lim}\ts \A_m, \qquad m\to\infty
$$
and it turns out that there exists a natural embedding
$
\U(\g(\infty))\to \A.
$
Finally, we show that an irreducible highest weight
representation $L_\lambda$ of the Lie algebra $\g(\infty)$ 
can be canonically extended to the algebra $\A$ provided the highest
weight $\lambda$ satisfies a stability condition described below; see
\thetag{0.5}. 

In contrast to the enveloping algebra $\U(\g(\infty))$, the algebra
$\A$ has a large center $\Z$, which coincides with the subalgebra 
$\A_0$ or $\A_{-1}$ for the series $C,D$ or $B$, respectively,
while the center of $\U(\g(\infty))$ is trivial. Thus, $\A$ 
is more like the enveloping algebras $\U(\g(n))$ than
$\U(\g(\infty))$. This idea is developed in the papers [O1],
[O2] which deal with the type $A$ case. In particular, it is explained
there how one can realize elements of $\A$ as left invariant
differential operators on a certain infinite-dimensional classical
group; then elements of the center $\Z$ become biinvariant (or Laplace)
operators.

The {\it stability condition\/} on $\lambda$ mentioned above takes
the following form: assume (to simplify the discussion) that 
$\g(\infty)$ is of type $C$ or $D$ and choose as a Cartan subalgebra
the set of doubly infinite diagonal matrices of the form
$$
a=\text{diag}\ts(\dots, 0,0,\dots,a_2,a_1,-a_1,-a_2,\dots,0,0,\dots);
$$
then
$$
\langle\lambda,a\rangle=\sum_{i=1}^\infty \lambda_i a_i,
$$
where
$$
\lambda_i=\lambda_{i+1}=\dots=c,\qquad i\gg 1,\tag 0.5
$$
which means that $\lambda$ has only a finite number of nonzero labels
on the (infinite) Dynkin diagram corresponding to the algebra
$\g(\infty)$. In particular, among the modules $L_\lambda$ satisfying
the stability condition there are analogs of finite-dimensional
modules. It is worth mentioning that the stable value $c$ of the
coordinates of $\lambda$ is involved as a parameter in the
centralizer construction. 

This paper is organized as follows. In Sections 1 and 2 
we generalize slightly
the construction of the algebra $\A\supset \U(\frak{gl}(\infty))$ from
[O2]. Here we 
make it depend upon a parameter. Although the
arguments do not differ substantially from those of [O2] we
present them in detail, to facilitate the reader's
understanding of the more involved proof in the case of type $B,C,D$
algebras. 

Section 3 is parallel to Section 1. 
Here we introduce a category $\Omega(c)$ of
$\g(\infty)$-modules, where $\g$ is of type $B, C$ or $D$, and $c$ is a
parameter. The category $\Omega(c)$ is similar to the category 
$\Cal O$ of Bernstein--Gelfand--Gelfand. We also introduce the
algebra $\Z$ of `virtual Laplace operators' (the future center of
$\A$) which operates in modules belonging to $\Omega(c)$.

The main results are concentrated in Section 4. We define the algebra
$\A$ (which depends on $c$), then show that it extends the action of
$\U(\g(\infty))$ in modules from the category $\Omega(c)$, and prove
the Main Theorem (Theorem 4.17)
about the structure of $\A$.

The results make it possible to study
representations of the centralizer algebras $\A_m(n)$
by using the techniques developed for twisted Yangians [M4].

The final Section 5 contains a few comments and open questions.

To conclude, let us mention some papers related to our subject.
A detailed exposition of the structure theory of the Yangian
$\Y(\frak{gl}(m))$ and the twisted Yangians is given in
[MNO]; this is our basic reference.
Representations of $\Y(\frak{gl}(m))$ are studied by 
Cherednik [C2],
Kirillov and Reshetikhin [KR], Molev [M1],
Nazarov and Tarasov [NT1], [NT2]. 
Applications of the twisted
Yangians to explicit constructions of central elements in classical
enveloping algebras and to Capelli identities are given in [M3] and [MN].
Families of Bethe-type
commutative subalgebras in $\Y^{\pm}(M)$ were constructed in [NO].

A remark is necessary in regard to our formula numbering: we enumerate 
the formulas in any subsection independently and use triple numbering when
referring them back in other subsections.

\bigskip
\bigskip
\noindent
{\bf 1. Highest weight modules and virtual Laplace
operators for the Lie algebra $\gl(\infty)$}
\bigskip

In this section we introduce a family of categories $\Omega(c)$,\ $c\in\C$
of modules over the Lie algebra 
$\gl(\infty)$. This category is
similar 
to
the well-known  category $\Cal O$ of Bernstein--Gelfand--Gelfand. We also
construct the algebras of `virtual Laplace operators' for this  Lie
algebra.

\bigskip
\noindent {\bf 1.1.} We denote by $\gln$ the Lie algebra of
complex $n\times n$-matrices. For any $n>1$ one has the natural inclusion
$\gl(n-1)\subset\gln$, where the subalgebra $\gl(n-1)$ is spanned by the
matrix units $E_{ij}$ with $1\leq i,j\leq n-1$. 
We denote by $\gl(\infty)$ 
the corresponding inductive limit of the Lie algebras $\gln$ as 
$n\to\infty$:
$$
\gl(\infty)=\underset {n=1}\to{\overset {\infty}\to {\bigcup}}\gln.
$$
In other words, $\gl(\infty)$ is the Lie algebra of all complex matrices
$A=(A_{ij})$, where $i$ and $j$ run over the set $\{1,2,\dots\}$, such that
the
number of nonzero entries $A_{ij}$ is finite. The matrix units $E_{ij},\ 
1\leq i,j<\infty$ form a basis of $\gl(\infty)$.

We let $\h$ denote the Cartan subalgebra of diagonal matrices and $\n_+$
(respectively, $\n_-$) the subalgebra of upper 
(respectively, lower) triangular matrices
in $\gl(\infty)$. One has the triangular decomposition
$$
\gl(\infty)=\n_-\oplus\h\oplus\n_+.
$$ 
For a linear functional $\l\in\h^*$ set $\l_i=\l(E_{ii})$.
We shall identify $\l$ with the sequence $(\l_1, \l_2,\dots)$.

\bigskip
\noindent {\bf 1.2.} A $\gl(\infty)$-module $V$ is said to be
{\it highest weight} if it has 
a cyclic vector $v$ satisfying
$
\n_+v=\{0\}
$
and there exists
$\l\in\h^*$ such that $hv=\l(h)v$ for any $h\in\h$. The functional $\l$
is the {\it highest weight\/} of $V$, and
$v$ is the {\it highest weight vector\/} of $V$; it is
unique up to scalar multiples. Sometimes $V$ is also referred to
as the {\it module with the highest weight\/} $\l$.
The universal $\gl(\infty)$-module $M_{\l}$ with the
highest weight
$\l\in\h^*$ (the {\it Verma module}) may be defined as the quotient of
the universal enveloping algebra $\U(\gl(\infty))$ by the left ideal
generated by $\n_+$ and the elements $h-\l(h),\ h\in\h$.
Denote by $L_{\l}$ the unique nontrivial irreducible  quotient
of $M_{\l}$.

Let $n$ be a positive integer.
For a $\gl(\infty)$-module $V$ and a functional 
$\mu=(\mu_1,\mu_2,\dots)\in\h^*$ set
$$
\align
V_n(\mu)=\{v\in V|\ &E_{ii}v=\mu_iv\quad
\text{for}\quad i>n\quad\text{and}\\
&E_{ij}v=0\quad\text{for}\quad 1\leq i<j;\ \ j>n\}.
\endalign
$$
Clearly, the subspaces $V_n(\mu)$ form an ascending chain of subspaces
in $V$. Set
$$
V_{\infty}(\mu)=\underset {n=1}\to{\overset {\infty}\to {\bigcup}}V_n(\mu)
$$
and note that $V_{\infty}(\mu)$ is a submodule of $V$ since $V_n(\mu)$
is $\gl(n)$-invariant for every $n$. We shall write $\mu\sim\mu'$ if
$\mu_i=\mu'_i$ for $i\gg 0$. Then $\mu\sim\mu'$ implies $V_{\infty}(\mu)=
V_{\infty}(\mu')$.

\bigskip
\noindent {\bf 1.3. Definition.} Fix a number $c\in\C$ and consider
the functional $\l^c:=(c,c,\dots)\in\h^*$. Then $\Omega(c)$ is defined to
be
the category of $\gl(\infty)$-modules $V$ satisfying the condition
$V_{\infty}(\l^c)=V$.

\bigskip
\proclaim {\bf 1.4. Proposition} Let $V$ be a $\gl(\infty)$-module
with the highest weight $\l$, where $\l\sim\l^c$. Then $V$ belongs to
$\Omega(c)$. In particular, $\Omega(c)$ contains the modules $M_{\l}$ and
$L_{\l}$ with $\l\sim\l^c$.
\endproclaim

\Proof Suppose that $\l_i=c$ for $i\geq n+1$. Then the highest weight
vector $v$ lies in the subspace $V_n(\l^c)$. Since $v$ is a cyclic vector
of $V$, the submodule $V_{\infty}(\l^c)$ coincides with $V$. 
$\square$

\bigskip
\noindent
{\bf 1.5. Remark.} In the category $\Omega(c)$, there is a distinguished
module,
namely, the one-dimensional module
$L_{(c,c,\ldots)}= L_{\l^c}$. For any complex numbers $c$ and $d$, the
mapping
$V\mapsto V\otimes L_{(d-c,d-c,\ldots)}$ establishes a category isomorphism
$\Omega(c)\to\Omega(d)$ such that
$L_{(c,c,\ldots)}\mapsto L_{(d,d,\ldots)}$. Thus, all the categories
$\Omega(c)$
are canonically isomorphic, and we could assume, with no real loss of
generality, that $c=0$. However, we prefer to keep the parameter $c$ in the
formulas below in order to emphasize a similarity with the case of Lie
algebras
$\oa(2\infty)$,
$\oa(2\infty+1)$, and $\spa(2\infty)$.

\bigskip
\noindent {\bf 1.6.} Let $\l\in\h^*$ be arbitrary.
For $n=1,2,\ldots$ denote by
$L_{\l}^{(n)}$ the cyclic $\gln$-submodule
of $L_{\l}$ generated by the highest
weight vector $v\in L_{\l}$. Set $\n_+^{(n)}=\n_+\cap\gln$. Then $v$
is annihilated by $\n_+^{(n)}$, so that $L_{\l}^{(n)}$ 
is a highest weight module over $\gln$.

\proclaim{\bf Proposition} The $\gln$-module $L_{\l}^{(n)}$
is irreducible for each $n$.
\endproclaim

\Proof Consider the decomposition
$
\n_+=\n_+^{(n)}\oplus\overline\n_+^{(n)},
$
where 
$$
\overline\n_+^{(n)}:=\text{\rm span of\ }\{E_{ij} \mid i<j;\ j\ge n+1\}.
$$
Note that $\overline\n_+^{(n)}$ is normalized by
$\n_-^{(n)}:=\n_-\cap\gln$,
that is, $[\n_-^{(n)},\overline\n_+^{(n)}]\subseteq\overline\n_+^{(n)}$. 
Since
$L_{\l}^{(n)}$ coincides with $\U(\n_-^{(n)})v$ and $v$ is annihilated by
$\overline\n_+^{(n)}\subset\n_+$, this implies that $L_{\l}^{(n)}\subset
L_{\l}$
is annihilated by $\overline\n_+^{(n)}$.

Now assume we are given a vector $w\in L_{\l}^{(n)}$ that is annihilated by
$\n_+^{(n)}$. Then, by the above argument, $w$ is annihilated by the whole
algebra $\n_+$. Since $L_{\l}$ is irreducible, $w$ must lie in $\C v$
which implies that $L_{\l}^{(n)}$ is irreducible. $\square$
\medskip

The proposition means that
$L_{\l}$ may be regarded as the inductive limit
$\underset{\longrightarrow}\to{\lim}\ts L_{\l}^{(n)}$
as $n\to\infty$ of
irreducible highest weight $\gln$-modules.

\bigskip
\noindent {\bf 1.7.} A functional $\l\in\h^*$ will be called a {\it
dominant
weight for\/} $\gl(\infty)$ if $\l_i-\l_{i+1}\in\ZZ_+$ for
$i=1,2,\ldots$\ts. 
This is equivalent to say that for any $n$ the restriction
of
$\l$ to $\h^{(n)}:=\h\cap\gln$ is a dominant weight for $\gln$ in the usual
sense.

It follows from Proposition 1.6 that any module $L_{\l}$ corresponding to a
dominant weight $\l$ is the inductive limit of irreducible
finite-dimensional
$\gln$-modules.

The dominant weights $\l\sim\l^0$ may be identified with Young
diagrams, and the corresponding $\gl(\infty)$-modules
$L_{\l}$ may be viewed as infinite-dimensional analogs of the
irreducible
polynomial $\gln$-modules. These modules $L_{\l}$ are closely related to
the
so-called tame irreducible representations of the group $U(\infty)$, see
Olshanski [O2].

Let us associate to $\gl(\infty)$ the semi-infinite Dynkin diagram of type
$A$.
Then there is a bijective correspondence between the dominant weights
$\l\sim\l^c$ and the labelings of the nodes of that diagram by
nonnegative integers where only a finite number of them is nonzero:
\medskip
$$
\overset{\dsize\ \ \l_1-\l_2\quad\ \ \ \l_2-\l_3
\quad\quad\quad\quad\quad\ \ 
\l_n-\l_{n+1}\ \ \quad\quad\ \ \ \ \ }\to
{{\ssize\bigcirc\ }\text{------------}{\ \ssize\bigcirc\ }
\text{---------}\cdots\text{---------}
{\ \ssize\bigcirc\ }\text{---------}\cdots}.
$$

\bigskip
\noindent 
{\bf 1.8.} Here we recall the well-known construction of the
Harish--Chandra isomorphism (see, e.g., Dixmier [D], Section 7.4).

Let $\gti$ be a reductive Lie algebra and $\hti$ a Cartan subalgebra of
$\gti$. Choose a system of positive roots for $(\gti,\hti)$ and consider
the
corresponding triangular decomposition
$\gti=\nti_-\oplus\hti\oplus\nti_+$ where $\nti_-$ and $\nti_+$ are spanned
by
the negative and positive root vectors, respectively. Then one has a
decomposition of the universal enveloping algebra
$$
\U(\gti)=(\nti_-\U(\gti)+\U(\gti)\nti_+)\oplus\U(\hti).\tag1
$$ 
Denote by $\U(\gti)^{\hti}$ the centralizer of $\hti$ in $\U(\gti)$. Then
the
projection onto the second component in (1) gives an algebra homomorphism
$\omega: \U(\gti)^{\hti}\to\U(\hti)$ called the {\it Harish--Chandra
homomorphism\/}.

Further, let $W$ be the Weyl group of $(\gti,\hti)$ and let
$\rho\in\hti^*$ be the half--sum of positive roots.
The following affine transformations of
$\hti^*$
$$
\l\mapsto w'(\l):=w(\l+\rho)-\rho,\qquad \l\in\hti^*,\quad w\in W\tag2
$$ 
form a finite group $W'$ isomorphic to $W$. The group $W'$ operates on the
algebra $\P(\hti^*)$ of polynomial functions on $\hti^*$,
which can be identified with $\U(\hti)$.

Let $\Z(\gti)$ denote the center of $\U(\gti)$. Since $\Z(\gti)\subset
\U(\gti)^{\hti}$, we may restrict $\omega$ to $\Z(\gti)$, and it turns out
that
$\omega: \;\Z(\gti)\to\U(\hti)=\P(\hti^*)$ is an isomorphism of the
algebra
$\Z(\gti)$ onto the algebra $\P(\hti^*)^{W'}$ of
$W'$-invariants in $\P(\hti^*)$. This isomorphism
$$
\omega:\Z(\gti)\to \P(\hti^*)^{W'}\tag3
$$ is called the {\it Harish--Chandra isomorphism\/}.

The Harish--Chandra isomorphism (3) has a natural interpretation in terms
of
eigenvalues of central elements.  Namely, let $a$ be an
arbitrary element of $\Z(\gti)$ and let
$f_a=\omega(a)\in \P(\hti^*)$ be the corresponding $W'$-invariant
polynomial.
Then $a$ operates in any highest weight $\gti$-module $V$ as the scalar
operator
$f_a(\l)\cdot 1$, where $\l$ stands for the highest weight of $V$.

\bigskip
\noindent {\bf 1.9.} Now we apply the above construction to $\gti=\gln$.
Take
$\nti_{\pm}=\n_{\pm}\cap\gln$, $\hti=\h^{(n)}=\h\cap\gln$ and identify
$\P(\hti^*)=\U(\hti)$ with $\C[\l_1,\ldots,\l_n]$, where the coordinates
$\l_1,\ldots,\l_n$ in $\h^{(n)}$ correspond to the basis
$E_{11},\ldots, E_{nn}$. Then
$$
\rho=\left(\frac{n-1}2,\frac{n-3}2,\ldots,-\frac{n-1}2\right)
=\left(\frac{n+1}2,\ldots,\frac{n+1}2\right) +(-1,-2,\ldots,-n).
$$ 
The Weyl group coincides with the symmetric group of degree $n$ permuting
the
coordinates $\l_1,\ldots,\l_n$. Since the vector
${\dsize\left(\frac{n+1}2,\ldots,\frac{n+1}2\right)}$ is invariant under
any
permutation of coordinates, we may replace in (1.8.2) the vector
$\rho$ by the vector $(-1,-2,\ldots,-n)$. This implies that a polynomial in
$\l_1,\ldots,\l_n$ is $W'$-invariant if and only if it is symmetric with
respect
to the new variables $\l_1-1,\l_2-2,\ldots,
\l_n-n$. The last property will be referred to
as the {\it shifted symmetry\/}, and we will denote by
$\Q(n)\subset\C[\l_1,\ldots,\l_n]$ the subalgebra of all shifted symmetric
polynomials. Thus, the Harish--Chandra isomorphism
for the Lie algebra $\gln$ is
an algebra isomorphism
$\omega:\Z(\gln)\to\Q(n)$. Note that this isomorphism preserves the
filtration: that of $\Z(\gln)$ is inherited from $\U(\gln)$ and that of
$\Q(n)$
is determined by the usual degree of polynomials.

\bigskip
\noindent {\bf 1.10.} Fix $c\in\C$ and consider morphisms
$$
\pi_{n,c}:\Q(n)\to\Q(n-1),\quad n=1,2,\ldots\tag1
$$ 
defined by the formula
$$ 
(\pi_{n,c}(f))(\l_1,\ldots,\l_{n-1})=f(\l_1,\ldots,\l_{n-1},c),\qquad
f\in\Q(n).
$$ 
Since $\pi_{n,c}$ preserves the filtration, we can define the projective
limit in the category of filtered algebras:
$$
\Qc=\underset{\longleftarrow}\to{\lim}\ts\Q(n), \qquad n\to\infty.
$$
In other words, an element
$f\in\Qc$ is a sequence $(f_n)$ such that 
$$
f_n\in\Q(n)\quad \text{for}\quad
n=1,2,\ldots,\qquad\pi_{n,c}(f_n)=f_{n-1}
\quad \text{for}\quad n\ge2
$$
and 
$$
\deg f:=\underset{n}\to{\sup}\ts\deg
f_n<\infty.
$$
We call $\Qc$ the {\it algebra of shifted symmetric functions}
(with parameter $c$); cf. [OO1].

Note that elements of $\Qc$ are
well-defined functions on the set of all
sequences $\l=(\l_1,\l_2,\ldots)\sim\l^c$.
Let us emphasize that the definition
of the algebra $\Qc$, contrary to that of $\Q(n)$, depends on the value of
the parameter $c$.

The construction of the algebra $\Qc$ of shifted symmetric functions is
quite
similar to that of the algebra $\Lambda$ of symmetric functions (see
Macdonald
[M]), and both algebras are indeed closely related to each other.
Namely, it is easily seen that $\Lambda$ is naturally isomorphic to the
graded
algebra
$\gr\Qc$ associated with the filtered algebra $\Qc$.

In the algebra $\Qc$, one can define analogs of power
sums, elementary symmetric functions, and complete symmetric functions:
$$
\align
 p_m(\l)&=\sum_{k=1}^{\infty}((\l_k-k)^m-(c-k)^m),\qquad m=1,2,\ldots,\\
1+\sum_{m=1}^{\infty}e_m(\l)t^m&=\prod_{k=1}^{\infty}
\frac{1+(\l_k-k)t}{1+(c-k)t},\\ 
1+\sum_{m=1}^{\infty}h_m(\l)t^m&=\prod_{k=1}^{\infty}
\frac{1-(c-k)t}{1-(\l_k-k)t},
\endalign
$$ 
where $\l=(\l_1,\l_2,\ldots)\sim\l^c$.

The interrelations between these functions are exactly the same as in the
case
of symmetric functions; see Macdonald [M], Chapter I, (2.6), (2.10),
(2.10${}'$). As for the algebra $\Lambda$, each of the families $\{p_m\}$,
$\{e_m\}$, $\{h_m\}$ can be taken as a system of algebraically independent
generators of the algebra $\Qc$. Thus we may write
$$
\Qc=\C[p_1,p_2,\ldots]=\C[e_1,e_2,\ldots]=\C[h_1,h_2,\ldots].
$$

It should be noted that there are other ways to define 
the shifted symmetric functions
$p_m$, $e_m$, and $h_m$ which can be
more suitable
for certain reasons; see [OO1].
However, the above `naive' definition will be sufficient
for the purposes of this
article.

\bigskip
\noindent
{\bf 1.11.} Consider the morphisms
$$
\pi_{n,c}: \Z(\gl(n))\to\Z(\gl(n-1)),\quad n=1,2,\ldots,\tag1
$$
which correspond to morphisms (1.10.1) under the identification
$\Z(\gln)\simeq\Q(n)$.

\medskip
\noindent
{\bf Definition.} The {\it algebra of virtual Laplace operators\/}, 
denoted as 
$\Z_c=\Z_c(\gl(\infty))$,  is defined as the projective limit
of the sequence of the filtered algebras $\Z(\gl(n))$: 
$$
\C\overset{\pi_{1,c}}\to{\longleftarrow}
{\Z}\bigl(\gl(1)\bigr)
\overset{\pi_{2,c}}\to{\longleftarrow}
{\Z}\bigl(\gl(2)\bigr)
\overset{\pi_{3,c}}\to{\longleftarrow}
{\Z}\bigl(\gl(3)\bigr)
\longleftarrow\ldots\ts.
$$
\medskip

Clearly, the algebra $\Z_c$ can be identified with the algebra $\Qc$ of 
shifted symmetric functions.

\bigskip
\bigskip
\noindent
{\bf 2. The centralizer construction for $\gl(\infty)$}
\bigskip

Here we introduce an
algebra $\A$
associated with the Lie algebra $\gl(\infty)$.
Then we prove that $\A$ is
isomorphic to the tensor product of the algebra of virtual Laplace
operators $\Z_c$, defined in Section 1, and the Yangian $\Y(\infty)$.
The latter algebra is the inductive limit of the Yangians
$\Y(n)=\Y(\gl(n))$ (see Definition 2.14).

\bigskip
\noindent
{\bf 2.1.} As in the previous section fix
$c\in\C$ and suppose that $m\in\{0,1,2,\dots\}$.
We shall assume that $n\geq m$.
Let us introduce the following notation:

$\g_m(n)$ is the subalgebra in $\gl(n)$ spanned by
the matrix units $E_{ij}$ subject to the condition $m+1\leq i,j\leq n$;

$\A(n)$ is the universal enveloping algebra of $\gl(n)$; 

$\A_m(n)$ is the centralizer of $\g_m(n)$ in $\A(n)$; in particular,
$\A_0(n)$ is the center of $\A(n)$;

$\G_m(n)$ is the subgroup in $\GL(n)$ consisting
of matrices stabilizing the basis vectors 
$e_i,\ 1\leq i\leq m$ in $\C^n$
and the space spanned by the remaining vectors;

$\I(n)$ is the left ideal in $\A(n)$ generated by the elements
$E_{in}-\delta_{in}c,\ i=1,\dots,n$;

$\J(n)$ is the right ideal in $\A(n)$ generated by the
elements
$E_{ni}-\delta_{ni}c,\ i=1,\dots,n$;

$\A(n)^0$ is the centralizer of $E_{nn}$
in $\A(n)$.

\bigskip
\proclaim
{\bf 2.2. Proposition} Let $\LL=\I(n)\cap\A(n)^0$. Then

{\rm (i)} $\LL=\J(n)\cap\A(n)^0$, so that $\LL$ is a two-sided ideal 
of the
algebra $\A(n)^0$;

{\rm (ii)} the following decomposition holds
$$
\A(n)^0=\LL\oplus\A(n-1).
$$
\endproclaim

\Proof  By the 
Poincar\'e--Birkhoff--Witt theorem, any element $a$ of the algebra
$\A(n)$ can be uniquely written as the sum of elements of
the form
$$
E_{n1}^{p_1}\cdots E_{n,n-1}^{p_{n-1}}x
E_{1n}^{q_1}\cdots E_{n-1,n}^{q_{n-1}}(E_{nn}-c)^r,\quad x\in\A(n-1).
\tag1
$$
Note that 
$$
[E_{nn},E_{ni}]=E_{ni},\quad [E_{nn},E_{in}]=-E_{in},
\quad 1\le i< n,
$$
and
$$
[E_{nn},x]=0,\quad x\in\A(n-1).
$$
This implies that any element of the form (1) is an eigenvector for 
$\text{ad}(E_{nn})$ with the eigenvalue
$$
(p_1+\cdots+p_{n-1})-(q_1+\cdots+q_{n-1}).\tag2
$$
Thus, $a$ belongs to $\A(n)^0$ if and only if for any component (1) 
in the decomposition of $a$ with $x\ne0$
the corresponding expression (2) vanishes. 

Now  let $a\in\A(n)^0$. Then $a$ is the sum of components of the form 
(1) with
$$
p_1+\cdots+p_{n-1}=q_1+\cdots+q_{n-1}.\tag3
$$
Clearly, any such component belongs to $\I(n)\cap\J(n)$ provided 
either the 
expression (3) or the number $r$ is nonzero. Further, there is at most one 
(nonzero) component for which both the sum (3) and the number $r$ vanish: 
such a component is an element $x\in\A(n-1)$. $\square$

\bigskip
\noindent
{\bf 2.3.} Let $\pi_{n,c}$ stand for the projection
 of the algebra
$\A(n)^0$ onto $\A(n-1)$ with the kernel $\LL$. By Proposition 2.2, 
$\pi_{n,c}$
is an algebra homomorphism (cf. this definition of $\pi_{n,c}$
with that of the Harish-Chandra homomorphism in 1.8).

\proclaim
{\bf Proposition} Let $m<n$. Then the restriction of $\pi_{n,c}$ to
$\A_m(n)$
defines an algebra homomorphism 
$$
\pi_{n,c}: \A_m(n)\to\A_m(n-1).
$$
\endproclaim

\Proof Indeed,
the ideal $\I(n)$ is clearly invariant under the adjoint action of
$\gl(n-1)$ and, in particular, under that of
the subalgebra $\g_m(n-1)$. Therefore the image $\pi_{n,c}(a)$ of any
element $a\in\A_m(n)$ lies in $\A_m(n-1)$. $\square$

\bigskip
\noindent
{\bf 2.4.} Consider the sequence of homomorphisms defined in
Proposition 2.3:
$$
\A_m(m)
\overset{\pi_{m+1,c}}\to{\longleftarrow}
\A_m(m+1)
\longleftarrow\cdots
\overset{\pi_{n,c}}\to{\longleftarrow}
\A_m(n)\longleftarrow\cdots\ts. \tag1
$$
It follows from the proof of Proposition 2.3 that these morphisms
preserve the canonical filtrations of the universal enveloping algebras.

\medskip
\noindent
{\bf Definition.} The algebra $\A_m$ is defined to be the
projective limit of the sequence (1) of the algebras $\A_m(n)$, 
where the limit
is taken in the category of filtered associative algebras. 
(Let us emphasize that this
construction depends on the parameter $c$, which is omitted only to 
simplify the notation.)
\medskip

In other words, an element of the algebra $\A_m$ is a sequence
of the form
$a=(a_m,a_{m+1},\dots,a_n,\dots)$ where $a_n\in\A_m(n)$,
$\pi_{n,c}(a_n)=a_{n-1}$ for $n>m$, and
$$
\deg a:=\sup\Sb n\ge m\endSb \deg a_n<\infty,\tag2
$$
where
$\deg a_n$ denotes the degree of $a_n$ in  the universal
enveloping algebra $\A(n)$.

\bigskip
\noindent
{\bf 2.5.} Note that the homomorphisms $\pi_{n,c}$ are compatible with the
natural embeddings $\A_m(n)\hra\A_{m+1}(n)$, that is, the following
diagram is commutative:
$$
\CD
\A_m(n) @>>> \A_{m+1}(n)\\
@V\pi_{n,c}VV     @VV\pi_{n,c}V\\
\A_m(n-1) @>>> \A_{m+1}(n-1).
\endCD
$$
Therefore we can define an embedding $\A_m\hra\A_{m+1}$ as follows:
$$
(a_m,a_{m+1},a_{m+2},\dots)\mapsto(a_{m+1},a_{m+2},\dots).
$$
\medskip
\noindent
{\bf Definition.} The algebra $\A$ is the inductive limit of
the associative filtered algebras:
$
\A=\underset m\to{\bigcup}\ts\A_m.
$

The canonical embedding of the universal enveloping algebra
$\A(\infty)$ into $\A$ is defined as follows: to an element
$x\in\A(\infty)$ we associate the sequence $(x,x,\dots)$.
This embedding is well defined since $x$ belongs to some
$\A(m)$ and hence to $\A_m(n)$ for all $n>m$.

\bigskip
\proclaim
{\bf 2.6. Proposition} The center of the algebra $\A$ coincides with
$\A_0$.
\endproclaim

\Proof It is clear that $\A_0$ is contained in the center.
Conversely, suppose that $a=(a_m,a_{m+1},\dots)\in\A$ 
is a central element.
Then it  commutes with $\gl(\infty)\subset\A$
and so $a_n$ lies in the center of
$\A(n)$ for all $n\geq m$ which implies that $a\in\A_0$. $\square$
\medskip

Note that the morphism $\pi_{n,c}: \A_0(n)\to\A_0(n-1)$
defined in Proposition 2.3, coincides with (1.11.1).
Therefore, $\A_0$ is
naturally isomorphic to the algebra $\Z_c$ of virtual Laplace operators.

\medskip

In the next few subsections we describe the structure of 
the graded algebras $\gr\ts\A_m$ and $\gr\ts\A$ with respect to
the filtration given by (2.4.2).

\bigskip
\noindent
{\bf 2.7.} We shall use the following notation: 

$\P(n)$ is the symmetric
algebra of the space
$\gl(n)$; 

$\P_m(n)$ is the subalgebra of the elements
of $\P(n)$, which are invariant under the adjoint action of the
subalgebra $\g_m(n)$; 

$\I'(n)$ is the ideal in $\P(n)$ generated by the
elements $E_{in},\ i=1,\dots,n$;

$\J'(n)$ is the ideal in $\P(n)$ generated by the
elements $E_{ni},\ i=1,\dots,n$.
\medskip

We can repeat the constructions from 2.1--2.6, replacing 
$\A_m(n)$ by $\P_m(n)$, $\I(n)$ by $\I'(n)$, and $\J(n)$ by $\J'(n)$.
As a result we obtain a commutative graded algebra $\P_m$ which is
the projective limit of the commutative graded algebras $\P_m(n)$. Then we
can define the algebra $\P$ as the inductive limit of the algebras
$\P_m$. 

Note that the algebra $\P$ contains 
the symmetric algebra of the Lie algebra
$\gl(\infty)$, and the algebra $\P_m$
coincides with the subalgebra of $\g_m(\infty)$-invariants in $\P$.

Here the role of the algebra $\Qc$ is played by 
the algebra of
symmetric functions. In more detail,
$\P_0(n)$ coincides with 
the algebra of 
invariants of $n\times n$ matrices (with respect to the adjoint action) 
and may 
be identified with the algebra of symmetric polynomials in $n$ variables 
$x_1,\ldots,x_n$. Then $\P_0$ may be identified with the 
algebra of symmetric functions in infinitely many variables
$x_1,x_2,\ldots$\ts.

\bigskip
\proclaim
{\bf 2.8. Proposition} The canonical isomorphisms $\gr\A(n)\to\P(n)$
induce isomorphisms $\gr\A_m\to\P_m$ and $\gr\A\to\P$.
\endproclaim

\Proof Let $\left(\A^k(n)\right)$, $k=0,1,2,\ldots$ stand for the
canonical
filtration of the universal enveloping algebra $\A(n)$, and let
$\P(n)=\oplus\ts\P^k(n)$ stand for the canonical gradation of the symmetric
algebra $\P(n)$. Set
$$
\A_m^k(n)=\A_m(n)\cap\A^k(n),\qquad
\P_m^k(n)=\P_m(n)\cap\P^k(n).
$$ 
The canonical isomorphism
$$
\A^k(n)/\A^{k-1}(n)\to\P^k(n)\tag1
$$ commutes with the adjoint action of $\gl(n)$
and hence with that of $\g_m(n)$. Since all the spaces in (1)
are
semisimple $\g_m(n)$-modules, we obtain from (1) the isomorphisms 
$$
\A_m^k(n)/\A_m^{k-1}\to\P_m^k(n),
$$ 
which allow us to identify ${\gr}\A_m(n)$ with $\P_m(n)$. 

Further, let $k$ be fixed. Then for any $n>m$ the diagram
$$
\CD
\A_m^k(n) @>>> \P_m^k(n)\\ @VVV @VVV\\
\A_m^k(n-1) @>>> \P_m^k(n-1)
\endCD
$$ 
is commutative, which follows immediately
from the definition of the ideals
$\I(n)$ and $\I'(n)$. This provides us with isomorphisms
${\gr}\A_m\to\P_m$ for any $m$, hence, with an isomorphism
${\gr}\A\to\P$. $\square$

\bigskip
\noindent
{\bf 2.9.} Let us identify $\gl(n)$ with its dual space using the 
bilinear form
$(X,Y)\mapsto\tr(XY^{\sigma})$ on $\gl(n)$, where $Y\mapsto Y^{\sigma}$ 
denotes
the standard matrix transposition:
$(E_{ij})^{\sigma}=E_{ji}$. Then we may identify the algebra $\P(n)$ with
the
algebra of polynomial functions $\phi(x)$ on $\gl(n)$. Under this
identification a matrix unit $E_{ij}\in\gl(n)$ becomes the function
$x\mapsto
x_{ij}$, where $x$ stands for a matrix of size
$n\times n$. (This is the reason why we have preferred to use the form
$\tr(XY^{\sigma})$ instead of the invariant form $\tr(XY)$).
 
\proclaim
{\bf Proposition} The algebra $\P_m(n)$ is generated by the
polynomials
$$
p_n^{(M)}(x)=\tr(x^M)  \tag1
$$
and
$$
p_{ij|n}^{(M)}(x)=(x^M)_{ij}, \tag2
$$
where $i,j\leq m$ and $M=1,2,\dots$\ts.
\endproclaim

\Proof For $m=0$ there are no elements of type (2) and
the claim is well known (see, e.g., [D], Section 7.3),
so we shall assume $m\geq 1$.
Write $x$ as a block matrix ${\dsize\binom{a\ \ b}{c\ \ d}}$ according
to the decomposition $n=m+(n-m)$. It suffices to prove that any
polynomial $\phi(x)$ satisfying the $\G_m(n)$-invariance condition
$$
\phi(x)\equiv\phi(a,b,c,d)=\phi(a,bg^{-1},gc,gdg^{-1}),
$$
where $g$ is an arbitrary invertible matrix of the same shape as $d$,
is a polynomial function of the invariants $\tr(x^M)$ and $(x^M)_{ij}$
with $i,j\leq m$.

Let us note first that these invariants can be replaced by the invariants
of 
the form $\tr(d^M),\ M\geq 1$,\ \ $(bd^{M-2}c)_{ij},\ M\geq 2$, and
$a_{ij}$.
Indeed, we have $x_{ij}=a_{ij}$, and for $M\geq 2$ the invariant
$(x^M)_{ij}-(bd^{M-2}c)_{ij}$ can be expressed in terms of the
invariants $a_{kl}$ and $(bd^{N-2}c)_{kl}$ with $N<M$. Hence,
instead of $(x^M)_{ij}$ we can take $(bd^{M-2}c)_{ij}$ for $M\geq 2$.
Moreover, $\tr(x^M)-\tr(d^M)$ can be expressed in terms of the
invariants of the form $a_{kl}$, $\tr(d^N)$ with $N<M$, and
$(bd^{N-2}c)_{kl}$ with $N\leq M$. Hence,
instead of $\tr(x^M)$ we can take $\tr(d^M)$.

To simplify the notation, we shall assume below that $m=1$ (the
generalization to the case $m>1$ will be obvious). We can ignore the
element $a$, so we have to show that any polynomial invariant
$\phi(b,c,d)$, where $c$ is an element of the $\GL(n-1)$-module $V$
of column vectors of the length $n-1$, $b$ is an element of the dual 
module $V^*$ of row vectors, and $d\in V\otimes V^*$, can be
expressed in terms of $\tr(d^M)$ and $bd^{M-2}c$.
Any such invariant can be decomposed into a sum of expressions of the 
form
$$
\psi(\undersetbrace p\to{b,\dots,b},
\undersetbrace q\to{c,\dots,c},\undersetbrace r\to{d,\dots,d}),\tag3
$$
where $\psi$ is a multilinear invariant. In its turn, $\psi$ is determined
by a multilinear invariant of the form
$$
\chi(b_1,\dots,b_p,c_1,\dots,c_q,u_1,\dots,u_r,v_1,\dots,v_r),\tag4
$$
where $b_1,\dots,b_p,u_1,\dots,u_r\in V^*$ and
$c_1,\dots,c_q,v_1,\dots,v_r\in V$. 

The multilinear invariants for the general linear group are described by
the classical invariant theory. Let us recall some basic facts of this
theory [W]. 
First of all, nonzero invariants exist only when the number of  
vector arguments is equal to the number of covector arguments. So in (4) 
we have $p=q$.  Next, any invariant can be uniquely written as a
polynomial in `elementary invariants' that come from the canonical 
pairing $V\otimes V^*\to\C$. In our notation, this means that $\chi$ 
may be represented as a linear combination of monomials in bilinear 
invariants
$b_ic_j,\ b_iv_j,\ u_ic_j,\ u_iv_j$ such that in
each of the monomials each letter appears exactly once. Since the letters
$u$ and $v$ are related with the block $d$, we must keep an eye on the
position of these letters with equal numbers. In each monomial we can
meet, first, closed chains of the form 
$$
(u_{k_1}v_{k_2})(u_{k_2}v_{k_3})\dots(u_{k_M}v_{k_1}),\tag5
$$
and, second, open chains of the form
$$
(b_kv_{m_1})(u_{m_1}v_{m_2})\dots(u_{m_{M-1}}c_l).\tag6
$$
Looking again at formula (3), we see that  under the passage 
$\chi\mapsto\psi\mapsto\phi$, each chain (5) gives rise to
the invariant $\tr(d^M)$, and each chain (6) corresponds
to an invariant of the form $bd^{M-2}c$. $\square$

\bigskip
\noindent
{\bf 2.10.} Fix numbers $K$ and $m$ such that $K\in\{1,2,\dots\}$ and
$m\in\{0,1,\dots\}$.
Assume that the indices $i,j,M$ satisfy the conditions $i,j\leq m$
and $1\leq M\leq K$.

\proclaim
{\bf Proposition} For a sufficiently large $n$ the elements
$p_{ij|n}^{(M)},\ p_n^{(M)}$ of the algebra $\P_m(n)$ with
the indices satisfying the above conditions are algebraically
independent.
\endproclaim

\Proof For any triple $(i,j,M)$ 
satisfying our assumptions we choose
a subset 
$$
\Omega_{ijM}\subset\{m+1, m+2,\dots\}
$$
of cardinality $M-1$ in such a way that all these subsets be disjoint.
Let $n$ be so large that all of them belong to $\{m+1, m+2,\dots,n-K\}$.
Let 
$$
y_i,\quad 1\leq i \leq K,\qquad\text{and}\qquad z_{ijM},\quad  
1\leq i,j\leq m,\quad  1\leq M\leq K,
$$
be complex parameters. Given $i,j,M$ we introduce a linear
operator $x_{ijM}$ in $\C^{n}$, which depends on the parameter $z_{ijM}$. 
Let $a_1<\dots<a_{M-1}$ be all the  elements of $\Omega_{ijM}$. Then
$x_{ijM}$ transforms the canonical basis vectors of $\C^n$ as follows:
$$
\gather
e_j\mapsto z_{ijM}e_{a_{M-1}},\\
e_{a_{M-1}}\mapsto e_{a_{M-2}},\ \ldots, e_{a_1}\mapsto e_i,\\
e_k\mapsto 0\quad \text{for}\quad  k\notin \{j\}\cup
\Omega_{ijM}. 
\endgather
$$
Next we define a linear operator $x$ in $\C^{n}$ depending on all 
the parameters,
\medskip
$$
xe_k=\cases \sum_{i,j,M}  x_{ijM}e_k,\quad &k\leq n-K,\\
y_{k-(n-K)}e_k,&k=n-K+1,\dots,n.\endcases\tag1
$$
\bigskip
\noindent
Regard $x$ as a $n\times n$ matrix. Then we have
$$
\gather
p_n^{(M)}(x)=y_1^M+\dots+y_K^M+\phi_M(\ldots,\  z_{klL},\ldots),\\
p_{ij|n}^{(M)}(x)=z_{ijM}+\psi_{ijM}(\ldots,\  z_{klL},\ldots),
\qquad L<M, 
\endgather
$$
where $\phi_M$ does not depend on $y_1,\dots,y_K$, and $\psi_{ijM}$
only depends on the parameters $z_{klL}$ with $L<M$. This implies that our
polynomials are algebraically independent even if 
they are restricted to the affine subspace
of matrices of the form (1). $\square$

\bigskip
\proclaim{\bf 2.11. Proposition} We have
$$
\gather
p_n^{(M)}\in p_{n-1}^{(M)}+\I'(n),\tag1\\
p_{ij|n}^{(M)}\in p_{ij|n-1}^{(M)}+\I'(n),\quad i,j\le n-1.\tag2
\endgather
$$
\endproclaim

\Proof By (2.9.1),
$$
p_n^{(M)}(x)=\tr(x^M)=\sum_{i_1,\ldots,i_M=1}^n
x_{i_1i_2}x_{i_2i_3}\cdots x_{i_Mi_1}.\tag3
$$
Here  $p_n^{(M)}$ is viewed as a polynomial function on the space of 
$n\times n$ matrices. 
Using the identification $x_{ij}\leftrightarrow E_{ij}$
(see 2.9)
we may rewrite (3) as follows
$$ 
p_n^{(M)}=\sum_{i_1,\ldots,i_M=1}^n E_{i_1i_2}E_{i_2i_3}\cdots
E_{i_Mi_1},\tag4 
$$
regarding $p_n^{(M)}$ as an element of $\P(n)$.

Now let us split the sum (4) into two parts. 
The first part will consist of all monomials with 
$i_1\leq n-1,\ldots,i_M\leq n-1$, 
and the second part will consist of all remaining monomials. Then the
first part of the sum clearly equals $p_{n-1}^{(M)}$. In the second part of
the sum each monomial contains a letter of the form $E_{jn}$ and so lies
in the ideal $\I'(n)$. This proves (1).

To verify (2), we write $p_{ij|n}^{(M)}$ as
$$ 
p_{ij|n}^{(M)}=\sum_{k_1,\ldots,k_{M-1}} E_{ik_1}E_{k_1k_2}\cdots
E_{k_{M-1}j},
$$
where the indices range over $\{1,\dots,n\}$
and repeat the previous argument. $\square$

\bigskip
\noindent
{\bf 2.12.} Given $M=1,2,\dots$, consider the sequence
$$
p^{(M)}:=(p_{n}^{(M)}),\qquad n\geq 1.
$$
We claim that $p^{(M)}$ is an element of $\P_0$.
Indeed, $p_{n}^{(M)}$ is a homogeneous element of $\P_0(n)$
of degree $M$, and by Proposition 2.11 the difference
$p_{n}^{(M)}-p_{n-1}^{(M)}$ lies in the ideal $\I'(n)$.
Similarly, fix $M,i,j,m$ such that
$
1\leq i,j\leq m.
$
Then the sequence
$$
p_{ij}^{(M)}:=(p_{ij|n}^{(M)}|n\geq m)
$$
is an element of $\P_m$.

\medskip
\proclaim{\bf Proposition} For any fixed $m\geq 0$
the elements $p^{(M)}$ and $p_{ij}^{(M)}$ with
$M=1,2,\dots$ and $1\leq i,j\leq m$, are algebraically
independent generators of the algebra $\P_m$.
\endproclaim

\Proof Let $p=(p_n|\ n\geq m)\in\P_m$. By Proposition 2.9, for any
$n\geq m$ the element $p_n\in \P_m(n)$ can be represented as a
polynomial $\phi_n$ in the variables $p_n^{(M)},\ p_{ij|n}^{(M)}$, where
$M\leq \deg p$ and $i,j\leq m$. 
By Proposition 2.10, $\phi_n$ does not 
depend on $n$ for sufficiently large $n$. Hence, the elements 
$p^{(M)}$ and $p_{ij}^{(M)}$ are generators. Their algebraic
independence is clear from Proposition 2.10. $\square$

\bigskip
\proclaim
{\bf 2.13. Corollary} The algebra $\P$ is isomorphic to the algebra
of polynomials in countably many variables $p^{(M)}$, $p_{ij}^{(M)}$,
where $M=1,2,\dots$ and $i,j\in\{1,2,\dots\}$. $\square$
\endproclaim

\medskip
Now we need some preparations to prove the main result of this
section, Theorem 2.18.
\bigskip
\noindent
{\bf 2.14.} 
Recall the definition of the {\it Yangian\/} $\Y(n)=\Y(\gl(n))$ 
(fore more details see
[MNO], Section~1).
It is a complex associative algebra with countably many generators
$$
t_{ij}^{(M)},\qquad 1\leq i,j\leq n,\qquad M=1,2,\dots\ts \tag 1
$$
and the quadratic defining
relations which can be written as follows.
Combine the generators
(1) into series
in $u^{-1}$:
$$
t_{ij}(u):=\delta_{ij}+\sum_{M=1}^{\infty}t_{ij}^{(M)}u^{-M},\qquad
1\leq i,j\leq n. \tag 2
$$
Then the defining relations take the form:
$$
[t_{ij}(u),t_{kl}(v)]={1\over u-v}(t_{kj}(u)t_{il}(v)-t_{kj}(v)t_{il}(u)).
\tag 3
$$
Next, the series (2) are combined into a single $n\times n$-matrix
$$
T(u)=\left(t_{ij}(u)\right)_{i,j=1}^n=\sum_{i,j=1}^n t_{ij}(u)\ot E_{ij}. 
\tag 4
$$
In terms of the $T$-matrix (4), the
relations (3)
can be written as a single `ternary relation' (0.2); see Introduction.

The {\it quantum determinant\/} $\qdet T(u)$
of the matrix $T(u)$ is a formal series in $u^{-1}$
with coefficients from $\Y(n)$ defined as follows: 
$$
\qdet T(u)=\sum_{p\in \Sym_n} \sgn(p)\ts t_{p(1),1}(u)\cdots
t_{p(n),n}(u-n+1), \tag 5
$$
where $\Sym_n$ is the group of permutations of the indices
$\{1,\dots,n\}$. 
The coefficients of the quantum
determinant $\qdet T(u)$ are algebraically independent generators of
the center of the algebra $\Y(n)$.

\bigskip
\noindent
{\bf 2.15.} Given $n$, let $E$ stand for the $\A(n)$-valued matrix
of order $n$ whose $(i,j)$-entry is
$E_{ij}\in\gl(n)\subset\A(n)$. (We use the symbol $E_{ij}$
to denote both the generators of $\gl(n)$ and the auxiliary matrix units).

Consider the mapping $\eta$ which takes the generators 
$t_{ij}^{(M)}$ of the
Yangian $\Y(n)$ to elements of $\A(n)$ as follows:
$$
\eta:T(u)\mapsto \left(1-\frac{E}{u}\right)^{-1},
$$
or, in more detail,
$$
\eta(t_{ij}^{(M)})=(E^M)_{ij}=
\sum_{k_1,\ldots,k_{M-1}} E_{ik_1}E_{k_1k_2}\cdots
E_{k_{M-1}j}\in\A(n),
$$
where $k_1,\ldots,k_{M-1}$ range over $\{1,\dots,n\}$.
It was proved in [MNO], Subsection~1.19 that
$\eta$ defines an algebra
homomorphism $\Y(n)\to \A(n)$. 

\bigskip
\proclaim{\bf 2.16. Proposition} For any $c\in\C$ the mapping
$$
\varphi_n: T(u)\mapsto {u+n-c\over u+n}
\left(1-{E\over u+n}\right)^{-1}\tag1
$$
defines an algebra homomorphism $\Y(n)\to \A(n)$.
\endproclaim

\Proof It was shown in [MNO], Proposition 1.12 that the following mappings
define automorphisms of the algebra $\Y(n)$:
shift of the formal parameter $u$ by a constant,
$
T(u)\mapsto T(u+a);
$
multiplication by a formal series,
$
T(u)\mapsto f(u)T(u),
$
where $f(u)=1+f_1u^{-1}+f_2u^{-2}+\cdots\ts$, $f_i\in\C$. Therefore, by
2.15 the mapping
$$
T(u)\mapsto f(u+a)\left(1-\frac {E}{u+a}\right)^{-1}
$$
defines an algebra homomorphism $\Y(n)\to \A(n)$. It remains to take
$f(u)=1-cu^{-1}$ and $a=n$. $\square$

\bigskip
\noindent
{\bf 2.17.} The defining relations (2.14.3) and the 
Poincar\'e--Birkhoff--Witt theorem for the Yangian (see [MNO],
Corollary 1.23) imply that for any $m\geq 1$ 
one has a
natural inclusion
$$
\Y(m)\hra\Y(m+1).\tag 1
$$
So, for any $m\leq n$ we can regard the Yangian $\Y(m)$ 
as a subalgebra in $\Y(n)$.

\proclaim
{\bf Proposition} The image of the restriction of the homomorphism
$\varphi_n$ to the subalgebra $\Y(m)$ is contained
in the centralizer $\A_m(n)$.
\endproclaim

\Proof It
follows from the defining relations (2.14.3) that
$$
[t_{kl}^{(1)}, t_{ij}(u)]=\delta_{il}t_{kj}(u)-
\delta_{kj}t_{il}(u).\tag 2
$$
In particular,
$$
[t_{kl}^{(1)}, t_{ij}(u)]=0 \tag 3
$$
for $i,j\leq m<k,l$.
On the other hand, it is easy to see from (2.16.1) that 
$$
\varphi_n(t_{kl}^{(1)})=E_{kl}-\delta_{kl}c. \tag 4
$$
Together with (3) this implies that
$
[E_{kl},\varphi_n(t_{ij}(u))]=0
$
provided that the indices $i,j,k,l$ satisfy the above restrictions.
By definition of $\A_m(n)$ this means that
$\varphi_n(t_{ij}^{(M)}) \in\A_m(n)$. $\square$

\bigskip
\noindent
{\bf  2.18.} Fix an arbitrary $m\geq 1$ and assume that
$n$ varies from $m$ to infinity. By Proposition 2.17, for any $n$
we have a homomorphism of the algebra $\Y(m)$ to
$\A_m(n)$, defined by $\varphi_n$.

\medskip
\proclaim
{\bf Theorem} For any fixed $m\geq 1$ the sequence
$(\varphi_n|n\geq m)$ defines an algebra embedding
$
\varphi: \Y(m)\hra\A_m.
$
Moreover, one has the isomorphism
$$
\A_m=\A_0\otimes\Y(m),\tag1
$$
where the Yangian is identified with its image under 
the embedding $\varphi$.
\endproclaim

\Proof We shall prove the theorem in several steps.
\smallskip
\noindent
{\sl Step} 1. To prove that the sequence of homomorphisms 
$(\varphi_n|n\geq m)$ defines an algebra homomorphism
$\varphi: \Y(m)\to\A_m$ we need to verify that
the following diagram is commutative:
$$
\CD
\Y(m) @= \Y(m) @= \cdots @= \Y(m)@=\cdots\\
@V \varphi_m VV     @V\varphi_{m+1}VV @. @V\varphi_{n}VV\\
\A_m(m) @<<\pi_{m+1,c}< \A_{m}(m+1)@<<<\cdots 
@<<\pi_{n,c}< \A_m(n)@<<<\cdots.
\endCD
$$
Denote the image of $t_{ij}(u)$ under the homomorphism (2.16.1)
by $\tau_{ij|n}(u)$ and
prove by induction on $M$ that for the coefficients of this series
one has (we use notation 2.1):
$$
\alignedat2
\text{(i)}\quad &\tau_{in|n}^{(M)}\in\I(n),
\quad &&1\leq i\leq n,\quad M\geq 1;\\
\text{(ii)}\quad &\tau_{ij|n}^{(M)}-\tau_{ij|n-1}^{(M)}\in\I(n),
\quad &&1\leq i,j\leq n-1,\quad M\geq 1.
\endalignedat
\tag2
$$
Consider the matrix $\T(u)=(\tau_{ij|n}(u))_{i,j=1}^n$. By
(2.16.1), 
$$
\T(u)(u+n-E)=u+n-c.
$$
Hence
$$
u\ts\T(u)=u+n-c+\T(u)(E-n).\tag3
$$
Write $\T(u)=\T^{(0)}+\T^{(1)}u^{-1}+\dots$\ts. Then (3) implies that
$\T^{(0)}=1$ and
$$
\align
\T^{(1)}&=E-c,\tag4\\
\T^{(M)}&=\T^{(M-1)}(E-n),\quad M\geq 2.\tag5
\endalign
$$
By (4) we have
$$
\tau_{in|n}^{(1)}=E_{in}-\delta_{in}c\in\I(n),\quad 1\leq i\leq n,
$$
and
$$
\tau_{ij|n}^{(1)}-\tau_{ij|n-1}^{(1)}=0,\quad 1\leq i,j\leq n-1.
$$ 
So, we have verified (2) for $M=1$. For $M>1$ we obtain from (5)
that
$$
\align
\tau_{in|n}^{(M)}={}&
\sum_{a=1}^n\tau_{ia|n}^{(M-1)}(E_{an}-\delta_{an}n)\\
{}={}&\sum_{a=1}^{n-1}\tau_{ia|n}^{(M-1)}E_{an}+
\tau_{in|n}^{(M-1)}(E_{nn}-c)+(c-n)\tau_{in|n}^{(M-1)}.
\endalign
$$
By the induction hypotheses, this expression lies in $\I(n)$, which
proves (i) in (2). Again using (5) we obtain for $1\leq i,j\leq n-1$ that
$$
\tau_{ij|n}^{(M)}-\tau_{ij|n-1}^{(M)}=
\sum_{a=1}^{n}\tau_{ia|n}^{(M-1)}(E_{aj}-\delta_{aj}n)
-\sum_{a=1}^{n-1}\tau_{ia|n-1}^{(M-1)}(E_{aj}-\delta_{aj}(n-1)),
$$
which can be rewritten as
$$
\sum_{a=1}^{n-1}(\tau_{ia|n}^{(M-1)}-\tau_{ia|n-1}^{(M-1)})E_{aj}
-(n-1)(\tau_{ij|n}^{(M-1)}-\tau_{ij|n-1}^{(M-1)})
+\tau_{in|n}^{(M-1)}E_{nj}-\tau_{ij|n}^{(M-1)}.
$$
Since $\varphi_n$ is an algebra homomorphism,
we obtain from (2.17.2) and (2.17.4) that
$$
\align
(\tau_{ia|n}^{(M-1)}-\tau_{ia|n-1}^{(M-1)})E_{aj}={}&
E_{aj}(\tau_{ia|n}^{(M-1)}-\tau_{ia|n-1}^{(M-1)})\\
{}+{}&\tau_{ij|n}^{(M-1)}-\tau_{ij|n-1}^{(M-1)}-
\delta_{ij}(\tau_{aa|n}^{(M-1)}-\tau_{aa|n-1}^{(M-1)}),
\endalign
$$
while
$$
\align
\tau_{in|n}^{(M-1)}E_{nj}-\tau_{ij|n}^{(M-1)}
={}&E_{nj}\tau_{in|n}^{(M-1)}+[\tau_{in|n}^{(M-1)},E_{nj}]
-\tau_{ij|n}^{(M-1)}\\
{}={}&E_{nj}\tau_{in|n}^{(M-1)}-\delta_{ij}\tau_{nn|n}^{(M-1)}.
\endalign
$$
Using (i) and the induction hypotheses, we complete the proof of (2).
In particular, we have proved that for $i,j\leq m$ the sequence
$(\tau_{ij|n}^{(M)}|\ n\geq m)$ is an element of the algebra $\A_m$,
and so, the homomorphism $\varphi$ is well defined.
\smallskip
\noindent
{\sl Step} 2. Here we verify that $\varphi$ is an embedding, that is,
its kernel is trivial.
Recall that $\A_m$ is a filtered algebra with $\gr\A_m=\P_m$;
see Proposition 2.8. Consider the highest order term of the
sequence $\tau_{ij}^{(M)}:=(\tau_{ij|n}^{(M)}|\ n\geq m)$.
By definition,
$$
\tau_{ij|n}(u)={u+n-c\over u+n}
\left(1-{E\over u+n}\right)_{ij}^{-1}.
$$
Hence 
$\tau_{ij|n}^{(M)}$ has the form
$$
\tau_{ij|n}^{(M)}=\bigl(E^M\bigr)_{ij}+\sum_{k\geq 1}a_k
\bigl(E^{M-k}\bigr)_{ij},
$$
where $a_k\in\C$. This proves that the image of $\tau_{ij|n}^{(M)}$ in
the $M$-th component of the graded algebra $\P(n)=\gr\A(n)$
is $p_{ij|n}^{(M)}$. So,
the highest order term of $\tau_{ij}^{(M)}$
coincides with $p_{ij}^{(M)}$.

Since the algebra $\gr\A_m=\P_m$ is commutative, we obtain from
Proposition 2.12 that the elements $\tau_{ij}^{(M)}$ of the algebra
$\A_m$ satisfy a Poincar\'e--Birkhoff--Witt-type condition: given an
arbitrary ordering 
of the elements $\tau_{ij}^{(M)}$,
any element of the subalgebra in $\A_m$ generated by
the $\tau_{ij}^{(M)}$ has a unique representation as a polynomial
in the $\tau_{ij}^{(M)}$. Using the Poincar\'e--Birkhoff--Witt theorem
for the Yangian (see [MNO], Corollary 1.23), we conclude that the
homomorphism $\varphi$ taking $t_{ij}^{(M)}$ to $\tau_{ij}^{(M)}$,
has zero kernel.
\smallskip
\noindent
{\sl Step} 3. It remains to prove the decomposition (1).
The arguments of Step 2 show that the graded algebra $\gr\Y(m)$
can be identified with the subalgebra
of $\gr\A_m=\P_m$ generated by the elements $p_{ij}^{(M)},
\ M=1,2,\dots$\ts. So, the required statement follows from
Proposition 2.12. $\square$

\bigskip
\noindent
{\bf 2.19.} By making use of the inclusion (2.17.1)
we define the algebra $\Y(\infty)$
as the corresponding inductive limit as $m\to\infty$.
The following corollary is
immediate from Theorem 2.18.

\proclaim
{\bf Corollary} One has the isomorphism
$
\A=\A_0\otimes\Y(\infty).\qquad \square
$
\endproclaim

\bigskip
\noindent
{\bf 2.20. Remark.} The Yangian $\Y(m)$ has a nontrivial center
generated by the coefficients of the quantum determinant $\qdet T(u)$;
see 2.14 and [MNO], Theorem 2.13.
Therefore, the center of the algebra $\A_m$ strictly contains $\A_0$.
However, the center of the algebra $\Y(\infty)$
is trivial, since by Proposition 2.6, the center of the algebra $\A$
coincides with $\A_0$. 

Note also that in contrast to the
algebra $\Y(m)$, the `infinite' Yangian
$\Y(\infty)$ apparently does not possess
a Hopf algebra structure.

\bigskip
\noindent
{\bf 2.21.} The following proposition shows that
the algebra $\A$ may be regarded as a `true' analog
of the universal enveloping algebra $\A(\infty)$ for the Lie algebra
$\gl(\infty)$.

\proclaim
{\bf Theorem} Any $\gl(\infty)$-module $V\in\Omega(c)$
has a natural $\A$-module structure.
\endproclaim

\Proof Let $v\in V$ and $a=(a_m,a_{m+1},\dots)\in\A$.
By (2.4.1), $a_{n+1}-a_n\in\I(n+1)$. Hence, if $n$ is large enough,
$a_nv$ does not depend on $n$ and this value can be taken as the
definition of the vector $av$. Let us verify that this definition
yields an action of $\A$. Indeed, if $a=(a_n)$ and $b=(b_n)$ are
elements of $\A$, then $abv$ is equal to $a_kb_lv$, where $k$ and $l$
are large enough. On the other hand, $(ab)v$ is $(a_nb_n)v=a_nb_nv$ for
a sufficiently large $n$, and so, $abv=(ab)v$.
This action of $\A$ clearly extends the action of its
subalgebra $\A(\infty)$. $\square$

\bigskip
\proclaim
{\bf 2.22. Proposition} Let $V$ be a $\gl(\infty)$-module
with the highest weight $\lambda$.
Then every element $a\in\A_0$ acts on $V$ as
the scalar operator $f_a(\lambda)\cdot 1$, where $f_a\in\Qc$ is the
shifted symmetric function which corresponds to $a$ under the
identification $\A_0=\Qc$ from {\rm 1.11}.
\endproclaim

\Proof Let $v$ be the highest weight vector of the module $V$ and
$a=(a_n)$ be an arbitrary element of $\A_0$.
The definitions of the Harish-Chandra isomorphism
(see 1.8) and the
isomorphism  $\A_0\to\Qc$ imply that
$a_nv=f_a(\l)v$ for a sufficiently large $n$. Hence $av=f_a(\l)v$.
Since $v$ is a cyclic vector and $a$ belongs to the center of 
the algebra $\A$, the same is true for all vectors of
the module $V$. $\square$

\bigskip
\noindent
{\bf 2.23. Remark.} The results of Sections 1 and 2 can be carried out
to the Lie algebra $\gl(2\infty+1)$ of all complex matrices $A=(A_{ij})$
where $i,j\in\ZZ$ and only a finite number of the $A_{ij}$ is nonzero.
In particular, here $\Omega(c)$ is replaced with the category
$\Omega(c_-,c_+)$ of modules which `stabilize' in both negative and
positive
directions (cf. 1.3, and 3.3 below). 
In contrast to the category $\Omega(c)$ for $\gl(\infty)$ (see Remark
1.5), the category $\Omega(c_-,c_+)$ for $\gl(2\infty+1)$ involves a
substantial continuous parameter, namely, the difference $c_--c_+$.

\bigskip
\bigskip
\noindent
{\bf 3. Highest weight modules and virtual Laplace
operators for the Lie algebras $\oa(2\infty)$, $\oa(2\infty+1)$ 
and $\spa(2\infty)$}
\bigskip

In this section we introduce analogs of the category $\Omega(c)$
for the infinite orthogonal and symplectic Lie algebras. 
Then we construct the algebras of `virtual Laplace operators' for 
these Lie algebras.

\bigskip
\noindent {\bf 3.1.} Let us introduce the Lie algebra $\gl(2\infty)$
of all complex matrices $A=(A_{ij})$ where $i$ and $j$ run through the
set $\ZZ\setminus\{0\}$ and the number of nonzero $A_{ij}$ is
finite. The standard matrix units $E_{ij},\ i,j\in\ZZ\setminus\{0\}$
form a basis of this algebra. Let us denote by $\g$ any of the following
Lie algebras:
$$
\align
\oa(2\infty)&=\{A\in\gl(2\infty)|\ A_{-j,-i}=-A_{ij}\},\\
\oa(2\infty+1)&=\{A\in\gl(2\infty+1)|\ A_{-j,-i}=-A_{ij}\},\\
\spa(2\infty)&=\{A\in\gl(2\infty)|\ A_{-j,-i}=
-\sgn i\ts\sgn j\cdot A_{ij}\}.
\endalign
$$
The Lie algebra $\g$ is spanned by the elements
$$
F_{ij}:=E_{ij}-\theta_{ij}E_{-j,-i},
$$
where $\theta_{ij}\equiv 1$ for the orthogonal 
Lie algebras, and $\theta_{ij}=\sgn i\ts\sgn j$ for 
the symplectic Lie algebra.

We let $\h$ denote the Cartan subalgebra of diagonal matrices,
$\n_+$ (respectively, $\n_-$) the subalgebra of upper 
(respectively, lower) triangular
matrices in $\g$, and $\frak b:=\h\oplus\n_+$ the Borel subalgebra
in $\g$. One has the triangular decomposition
$$
\g=\n_-\oplus\h\oplus\n_+.
$$
For a linear functional $\l\in\h^*$ set $\l_i=\l(F_{ii})$. Then we have
$\l_{-i}=-\l_i$, so $\l$ is uniquely determined by the sequence
$(\ldots,\l_{-2},\l_{-1})$.

\bigskip
\noindent
{\bf 3.2.} A $\g$-module $V$ is said to be
{\it highest weight} if it has 
a cyclic vector $v$ satisfying
$
\n_+v=\{0\}
$
and there exists
$\l\in\h^*$ such that $hv=\l(h)v$ for any $h\in\h$. The functional $\l$
is the {\it highest weight\/} of $V$, and
$v$ is the {\it highest weight vector\/} of $V$; it is
unique up to scalar multiples. Sometimes $V$ is also referred to
as the {\it module with the highest weight\/} $\l$.
The universal $\g$-module $M_{\l}$ with the
highest weight
$\l\in\h^*$ (the {\it Verma module}) is defined as the quotient of
the universal enveloping algebra $\U(\g)$ by the left ideal
generated by $\n_+$ and the elements $h-\l(h),\ h\in\h$.
Denote by $L_{\l}$ the unique nontrivial irreducible  quotient
of $M_{\l}$.

For $n=1,2,\dots$ set
$$
\align
\g(n)={}&\text{span of\ }\{F_{ij}|-n\leq i,j\leq n\}\subset\g,\\
\frak b_n={}&\{A\in\frak b|\ A_{ij}=0\quad\text{for}\quad
-n\leq i,j\leq n\}\subset\frak b.
\endalign
$$
Then $\g(n)$ is isomorphic to $\oa(2n),\ \oa(2n+1)$, 
or $\spa(2n)$, respectively, 
and $\g(n)+\frak b_n$ is a parabolic subalgebra of
$\g$.

For a $\g$-module $V$ and a functional $\mu\in\h^*$ define an ascending
chain of subspaces
$$
V_n(\mu)=\{v\in V|\ F_{ii}v=\mu_iv\quad\text{and}\quad
F_{ij}v=0,\ i<j,
\quad\text{for}\quad F_{ii}, F_{ij}\in\frak b_n\}.
$$
Set
$$
V_{\infty}(\mu)=\underset {n=1}\to{\overset {\infty}\to {\bigcup}}V_n(\mu)
$$
and note that $V_{\infty}(\mu)$ is a submodule of $V$, since $V_n(\mu)$
is $\g(n)$-invariant for every $n$. We shall write $\mu\sim\mu'$ if
$\mu_i=\mu'_i$ for $|i|\gg 0$. Then $\mu\sim\mu'$ implies $V_{\infty}(\mu)=
V_{\infty}(\mu')$.

\bigskip
\noindent
{\bf 3.3. Definition.} Fix a number $c\in\C$ and consider
the functional $\l^c\in\h^*$ such that $(\l^c)_i=c$ for $i=-1,-2,\dots$\ts.
We define the category $\Omega(c)$ as the category of those $\g$-modules
$V$ for which $V_{\infty}(\l^c)=V$.

\bigskip
\noindent
{\bf 3.4.} One easily checks that any $\g$-module
with the highest weight $\l\sim\l^c$ belongs to
$\Omega(c)$. In particular, $\Omega(c)$ contains the modules $M_{\l}$ and
$L_{\l}$ with $\l\sim\l^c$ (cf. Proposition 1.4).

\bigskip
\noindent
{\bf 3.5.} For an arbitrary $\l\in\h^*$ denote by $L_{\l}^{(n)}$
the cyclic $\g(n)$-submodule  of $L_{\l}$ generated by the highest weight
vector $v\in L_{\l}$. Then $L_{\l}^{(n)}$ is clearly a highest weight
module over $\g(n)$.
Exactly as in the $\gln$-case (see Proposition 1.6), one proves that
all $\g(n)$-modules $L_{\l}^{(n)}$ are irreducible, so that
$L_{\l}$ can be viewed as the inductive limit
$\underset{\longrightarrow}\to{\lim}\ts L_{\l}^{(n)}$ as $n\to\infty$.
 
Now we suppose that $\l\in\h^*$ satisfies the conditions
$$
\aligned
\text{(i)}\quad &\l\sim\l^c,\\
&\text{where}\quad c\in\tfrac12\ZZ_+
\quad\text{for}\quad \g=\oa(2\infty)\quad\text{or}\quad
\g=\oa(2\infty+1),\\
&\text{and}\quad c\in\ZZ_+
\quad\text{for}\quad \g=\spa(2\infty);\\
\text{(ii)}\quad &\l_{-(i+1)}-\l_{-i}\in\ZZ_+\quad\text{for}\quad
i=1,2,\dots,\\
&\text{and in addition}\quad \l_{-2}\geq |\l_{-1}|
\quad\text{for}\quad \g=\oa(2\infty),\\
&\text{and}\quad \l_{-1}\geq 0
\quad\text{for}\quad \g=\oa(2\infty+1)
\quad\text{and}\quad \g=\spa(2\infty).
\endaligned
$$
Then $L_{\l}^{(n)}$ is a finite-dimensional
irreducible module with the dominant
highest weight obtained by restricting $\l$ to $\h^{(n)}:=\h\cap\g(n)$.
So, the module $L_{\l}$ is the inductive limit of
finite-dimensional irreducible $\g(n)$-modules.

Let us associate a
semi-infinite Dynkin diagram with the Lie algebra $\g$. Just as
in the finite-dimensional case we can represent any $\l\in\h^*$
as a collection of labels attached to the nodes of this diagram
as follows:
$$
\overset{\dsize\overset{\dsize\ \ \quad\quad\quad\ \ \l_{-n}-\l_{-n+1}
\quad\quad\quad
\l_{-3}-\l_{-2} \quad\ \  \l_{-2}-\l_{-1} }
\to
{\cdots
\text{---------}{\ \ssize\bigcirc\ }\text{---------}
\cdots
\text{---------}
{\ \ssize\bigcirc\ }\text{--------------}{\ \ssize\bigcirc}}}
\to
{\quad\quad\quad\quad\quad\quad\quad\quad\quad\quad\quad\quad\ \ 
\overset{\big\vert\ \ \ \quad\quad\quad\quad}
\to{\ssize\bigcirc\ {\dsize\l_{-2}+\l_{-1}}}}
$$
in the case of $\g=\oa(2\infty)$;
$$
\overset{\dsize\ \ \quad\quad\quad\l_{-n}-\l_{-n+1}
\quad\quad\quad\ \ 
\l_{-3}-\l_{-2}\quad\quad\l_{-2}-\l_{-1}\quad\ \ \ 
2\l_{-1}\ \ }
\to
{\cdots
\text{---------}{\ \ssize\bigcirc\ }\text{---------}
\cdots
\text{---------}
{\ \ssize\bigcirc\ }\text{---------------}{\ \ssize\bigcirc}
\CD @= \endCD{\ssize\bigcirc}}
$$
in the case of $\g=\oa(2\infty+1)$;
$$
\overset{\dsize\ \ \quad\quad\quad\l_{-n}-\l_{-n+1}
\quad\quad\quad\ \ 
\l_{-3}-\l_{-2}\quad\quad\l_{-2}-\l_{-1}\quad\ \ \ 
\l_{-1}\ \ }
\to
{\cdots
\text{---------}{\ \ssize\bigcirc\ }\text{---------}
\cdots
\text{---------}
{\ \ssize\bigcirc\ }\text{---------------}
{\ \ssize\bigcirc}\CD @= \endCD
{\ssize\bigcirc}}
$$
in the case of $\g=\spa(2\infty)$.
Then the conditions (i) and (ii) express the fact that all the labels
are nonnegative integers and the number 
of nonzero labels is finite.

\bigskip
\noindent
{\bf 3.6.} Let $\rho\in\h^*$ denote the half sum of positive roots
of $\g$ relative to $\frak b$, or which is the same, $\rho|_{\h^{(n)}}$
is the half sum of positive roots of $\g(n)$ for any $n$. Then for
$i=1,2,\dots$ we have
$$
\rho_{-i}=-\rho_i=\cases i-1,&\text{for$\quad\g=\oa(2\infty),$}\\
i-\tfrac12,&\text{for$\quad\g=\oa(2\infty+1),$}\\
i,&\text{for$\quad\g=\spa(2\infty)$}.\endcases
$$
Denote by $\Z(n)$ the center of the universal enveloping algebra
$\U(\g(n))$. Using the Harish-Chandra isomorphism (see 1.8), we
identify $\Z(n)$ with the algebra $\M(n)$ of the polynomial
functions $f(\l_{-n},\dots,\l_{-1})$ on $(\h^{(n)})^*$ which are
invariant under the `shifted' action (1.8.2)
of the Weyl group $W$. More precisely,
set $l_i=\l_i+\rho_i$ and consider $f$ as a function in the
variables $l_{-n},\dots,l_{-1}$. Then $f$ must be invariant
with respect to all permutations of the variables and all
transformations $l_{-i}\mapsto\pm l_{-i}$, where in the
case of $\g(n)=\oa(2n)$ the number of `$-$' has to be even.
Note that both canonical filtrations of $\U(\g(n))$ and
$\C[\l_{-n},\dots,\l_{-1}]$ induce the same filtration
on $\Z(n)$.

\bigskip
\noindent
{\bf 3.7.} Let us fix $c\in\C$ and consider morphisms
$
\pi_{n,c}: \Z(n)\to\Z(n-1)
$
defined as follows:
$$
(\pi_{n,c}(f))(\l_{-n+1},\dots,\l_{-1})=f(c,\l_{-n+1},\dots,\l_{-1}),
\tag1
$$
where $f\in \Z(n)$ is identified with an element of $\M(n)$.
\medskip

\noindent
{\bf Definition.} The {\it algebra\/} $\Z_c=\Z_c(\g)$ {\it 
of virtual Laplace operators\/} is defined as the projective limit
of the sequence of the filtered algebras $\Z(n)$ as $n\to\infty$. 

\bigskip
\noindent
{\bf 3.8.} Let $\Mc$ denote the algebra
generated by the family of
functions
$p_m(\l),\ m=1,2,\dots$, where
$$
p_m(\l)=\sum_{i=-1,-2,\dots}(l_i^{\ts 2m}-(l_i^c)^{2m}),
$$
$l=\l+\rho$,\  $l^c=\l^c+\rho$, and $\l\sim \l^c$.

\proclaim
{\bf Proposition} The algebras $\Z_c$ and $\Mc$
are isomorphic to each other.
\endproclaim

\Proof Denote by $\Z^k(n)$ the subspace of elements of the algebra
$\Z(n)$ of degree $\leq k$ and by $\Z^k_c$ the space of all sequences
of the form $z=(z_1,z_2,\dots)$ such that $z_n\in\Z^k(n)$ and
$\pi_{n,c}(z_n)=z_{n-1}$. Then by Definition 3.7,
$$
\Z_c=\underset{k\geq 1}\to{\bigcup}\Z^k_c.
$$
It follows from the well-known description of the polynomial
$W$-invariants (see, e.g., \v{Z}elobenko [\v{Z}]) that the algebra
$\M(n)$ is generated by 
the polynomials
$$
p_{m|n}(\l)=\sum_{i=-n}^{-1}(l_i^{2m}-(l_i^c)^{2m}),
\quad m=1,2,\dots n,
$$
and in the case of $\oa(2n)$ also by the
polynomial
$p'_n(\l)=l_{-1}\cdots l_{-n}$.

Let $z=(z_n)\in \Z^k_c$. Identify $z_n$ with an element of $\M(n)$.
Then for $n>k$ this element can be
uniquely represented as a polynomial $P_n$ in
the $p_{m|n}$ with $2m\leq k$. (The element $p'_n$ 
in the case of $\oa(2n)$ is eliminated,
because $\deg p'_n>\deg z_n$). Since the polynomials
$p_{m|n}$ are algebraically independent and
$\pi_{n,c}(p_{m|n})=p_{m|n-1}$, we have
$\pi_{n,c}(P_n)=P_{n-1}$ for $n>k+1$.
This means that $z$ is
represented as a polynomial in the elements
$p_m:=(p_{m|n}|\ n\geq m)$. $\square$

\bigskip
\noindent
{\bf 3.9.} The algebra $\Mc$ can be regarded as an analog of the algebra of
shifted symmetric functions $\Lambda^*_c$ (see 1.10) for the case of
$B,C,D$ series. As with $\Lambda^*_c$ one can introduce analogs
of the elementary and complete symmetric functions by the formulas
$$
\align
1+\sum_{m=1}^{\infty}e_m(\l)t^m&=\prod_{k=1}^{\infty}
\frac{1+l^2_k t}{1+(l^c_k)^2 t},\\ 
1+\sum_{m=1}^{\infty}h_m(\l)t^m&=\prod_{k=1}^{\infty}
\frac{1-(l^c_k)^2 t}{1-l^2_k t}.
\endalign
$$ 
Each of the families $\{p_m\}$,
$\{e_m\}$, $\{h_m\}$ can be taken as a system of algebraically independent
generators of the algebra $\Mc$, so that
$$
\Mc=\C[p_1,p_2,\ldots]=\C[e_1,e_2,\ldots]=\C[h_1,h_2,\ldots].
$$

\bigskip
\bigskip
\noindent
{\bf 4. The centralizer construction for the
orthogonal and symplectic
Lie algebras}
\bigskip

Here we carry over the construction of Section 2 to each of
the Lie algebras 
$\g=\oa(2\infty),\oa(2\infty+1), \spa(2\infty)$.
We construct an algebra $\A=\A(\g)$
and prove that $\A$
is isomorphic to the tensor product of the algebra of virtual Laplace
operators $\Z_c$ and the algebra $\Y^{\pm}(\infty)$.
The latter algebra is the inductive limit of the twisted
Yangians $\Y^{\pm}(N)$
as $N\to\infty$; 
see 4.12 for the definition of $\Y^\pm(N)$. The key calculations
are contained in 4.14 and 4.16, and the result is stated in Theorem
4.17. The original proof of the theorem, as outlined in [O3],
Section 5, was simplified thanks to application of the automorphism
(4.14.7) of the twisted Yangian, which was introduced in [M4].

\bigskip
\noindent
{\bf 4.1.} Fix a number
$c\in\C$ and suppose that $m\in\{0,1,2,\dots\}$ in the case of
$\g=\oa(2\infty)$, $\spa(2\infty)$ and 
$m\in\{-1,0,1,\dots\}$ in the case of
$\g=\oa(2\infty+1)$. 
Assume that $n\geq m$.
We shall use the following notation:

$\g_m(n)$ is the subalgebra in $\g(n)$ spanned by
the elements $F_{ij}$ subject to the condition $m+1\leq |i|,|j|\leq n$;

$\A(n)$ is the universal enveloping algebra of $\g(n)$;

$\Z(n)$ is the center of $\A(n)$;

$\A_m(n)$ is the centralizer of $\g_m(n)$ in $\A(n)$;
in particular, $\Z(n)$ coincides with $\A_0(n)$ or $\A_{-1}(n)$
in the case of $\g=\oa(2\infty),\spa(2\infty)$ or
$\g=\oa(2\infty+1)$, respectively; 

$\G(n)$ is the classical Lie group corresponding to the 
Lie algebra $\g(n)$, namely $\G(n)=\SO(2n), \SO(2n+1)$, and $\SP(2n)$ for
$\g(n)=\oa(2n), \oa(2n+1)$, and $\spa(2n)$, respectively;

$\G_m(n)$ is the subgroup in $\G(n)$ consisting
of matrices stabilizing the basis vectors 
$e_i,\ |i|\leq m$ in the coordinate space
and the subspace spanned by the remaining basis vectors;

$\I(n)$ is the left ideal in $\A(n)$ generated by the elements
$F_{in}+\delta_{in}c,\ i=-n,\dots,n$;

$\J(n)$ is the right ideal in $\A(n)$ generated by the
elements
$F_{ni}+\delta_{ni}c,\ i=-n,\dots,n$;

$\A(n)^0$ is the centralizer of $F_{nn}$
in $\A(n)$ (note that $\A_m(n)\subset\A(n)^0,\ m<n$).

\bigskip
\proclaim
{\bf 4.2. Proposition} Let $\LL=\I(n)\cap\A(n)^0$. Then

{\rm (i)} $\LL=\J(n)\cap\A(n)^0$ and $\LL$ is a two-sided ideal of the
algebra $\A(n)^0$;

{\rm (ii)} one has the decomposition
$$
\A(n)^0=\LL\oplus\A(n-1).
\tag1
$$
\endproclaim

\Proof By the 
Poincar\'e--Birkhoff--Witt theorem, any element $a$ of the algebra
$\A(n)$ can be uniquely written as a sum of the components of
the form
$$
F_{n,-n+1}^{p_{-n+1}}\cdots F_{n,n-1}^{p_{n-1}}x
F_{n-1,n}^{q_{n-1}}\cdots F_{-n+1,n}^{q_{-n+1}}(F_{nn}+c)^r,\quad
x\in\A(n-1), \tag2
$$
in the orthogonal case; and
$$
F_{n,-n}^{p_{-n}}\cdots F_{n,n-1}^{p_{n-1}}x
F_{n-1,n}^{q_{n-1}}\cdots F_{-n,n}^{q_{-n}}(F_{nn}+c)^r,\quad x\in\A(n-1),
\tag3
$$
in the symplectic case.
The condition $a\in\A(n)^0$ means that $[a,F_{nn}]=0$. So, if
$x\ne 0$ then
$$
p_{-n+1}+\dots+p_{n-1}=q_{-n+1}+\dots+q_{n-1}, \tag4
$$
in the orthogonal case; and
$$
2p_{-n}+p_{-n+1}+\dots+p_{n-1}
=2q_{-n}+q_{-n+1}+\dots+q_{n-1}, \tag5
$$
in the symplectic case.
Hence the component (2) (respectively, (3)) belongs to $\I(n)\cap \J(n)$ 
if either $r\ne 0$ or the sum (4) (respectively, (5))
does not vanish. This proves (i). Furthermore, the component
(2) (respectively, (3)) belongs to $\A(n-1)$ 
if and only if $r=0$ and the sum (4) 
(respectively, (5)) vanishes, that is,
$p_i=q_i=0$ for all $i$, which implies (1). $\square$

\bigskip
\noindent
{\bf 4.3.} By Proposition 4.2, the projection $\pi_{n,c}$ of the algebra
$\A(n)^0$ onto $\A(n-1)$ with the kernel $\LL$ is an algebra homomorphism.

\proclaim
{\bf Proposition} Let $m<n$. The restriction of $\pi_{n,c}$ to
$\A_m(n)$
defines an algebra homomorphism
$$
\pi_{n,c}: \A_m(n)\to\A_m(n-1). 
$$
\endproclaim

\Proof Indeed,
the ideal $\I(n)$ is clearly invariant under the adjoint action of
$\g(n-1)$ and, in particular, under that of
the subalgebra $\g_m(n-1)$. Therefore the image $\pi_{n,c}(a)$ of any
element $a\in\A_m(n)$ lies in $\A_m(n-1)$. $\square$

\bigskip
\noindent
{\bf 4.4.} Consider the sequence of morphisms defined in
Proposition 4.3:
$$
\A_m(m)
\overset{\pi_{m+1,c}}\to{\longleftarrow}
\A_m(m+1)
\longleftarrow\cdots
\overset{\pi_{n,c}}\to{\longleftarrow}
\A_m(n)\longleftarrow\cdots\ts. \tag1
$$
It follows from the proof of Proposition 4.3 that these morphisms
preserve the canonical filtrations of universal enveloping algebras.

\medskip
\noindent
{\bf Definition.} The associative algebra $\A_m$ is defined as the
projective limit of the sequence (1) of the algebras $\A_m(n)$ in
the category of filtered algebras.
\medskip

An element of the algebra $\A_m$ is a sequence
$a=(a_m,a_{m+1},\dots,a_n,\dots)$, where $a_n\in\A_m(n)$ and
$\pi_{n,c}(a_n)=a_{n-1}$ for each $n$. Moreover,
$\text{sup}\ts \deg a_n <\infty$.
Here $\deg a_n$ denote the degree of the element $a_n$ in the universal
enveloping algebra $\A(n)$. For $a\in\A_m$ set
$$
\deg a:=\underset{n\geq m}\to{\text{sup}}\ts\deg a_n.
\tag2
$$

\bigskip
\noindent
{\bf 4.5.} Note that the morphisms $\pi_{n,c}$ are compatible with the
embeddings $\A_m(n)\hra\A_{m+1}(n)$ (cf. 2.5) and so,
we can define an embedding $\A_m\hra\A_{m+1}$ as follows:
$$
(a_m,a_{m+1},a_{m+2},\dots)\mapsto(a_{m+1},a_{m+2},\dots).
$$
\medskip
\noindent
{\bf Definition.} The algebra $\A$ is the inductive limit of
associative filtered algebras:
$
\A=\underset m\to{\bigcup}\ts\A_m.
$
\medskip

The canonical embedding of the universal enveloping algebra
$\A(\infty)$ into $\A$ is defined as follows: to an element
$x\in\A(\infty)$ we associate the sequence $(x,x,\dots)$.
This embedding is well defined since $x$ belongs to some
$\A(m)$ and hence to $\A_m(n)$ for all $n>m$.

\bigskip
\noindent
{\bf 4.6.} The same argument as in 2.6 proves that
the center of the algebra $\A$ coincides with
$\A_0$ in the case of $\g=\oa(2\infty),\spa(2\infty)$ 
and with $\A_{-1}$ in
the case of $\g=\oa(2\infty+1)$.

Note that the morphism $\pi_{n,c}: \A_0(n)\to\A_0(n-1)$ (respectively,
$\pi_{n,c}: \A_{-1}(n)\to\A_{-1}(n-1)$) 
defined in Proposition 4.3 coincides with the morphism (3.7.1).
Therefore, the algebra $\A_0$ (respectively, $\A_{-1}$) is
naturally isomorphic to the algebra of virtual Laplace operators 
$\Z_c=\Z_c(\g)$.

\medskip
In the next few subsections we describe the structure of 
the graded algebras $\gr\A_m$ and $\gr\A$ with respect to
the filtration given by (4.4.2).

\bigskip
\noindent
{\bf 4.7.} We shall use the following notation:

$\P(n)$ is the symmetric algebra of the space
$\g(n)$; 

$\P_m(n)$ is the subalgebra of the elements of $\P(n)$, which
are invariant under the adjoint action of the subalgebra $\g_m(n)$;

$\I'(n)$ is the ideal in $\P(n)$ generated by the elements $F_{in}$,
$i=-n,\dots,n$; 

$\J'(n)$ is the ideal in $\P(n)$
generated by the elements $F_{ni}$,
$i=-n,\dots,n$.
\medskip

We can repeat the constructions from 4.1--4.6, replacing
$\A_m(n)$ with $\P_m(n)$, $\I(n)$ with $\I'(n)$, and
$\J(n)$ with $\J'(n)$, thus obtaining a commutative graded algebra
$\P_m$ which is the projective limit of the commutative graded algebras
$\P_m(n)$. Then we can define the algebra $\P$ as the inductive limit
of the algebras $\P_m$.

Note that the algebra $\P$ contains 
the symmetric algebra of the Lie algebra
$\g(\infty)$, and the algebra $\P_m$
coincides with the subalgebra of $\g_m(\infty)$-invariants in $\P$.

\bigskip
\proclaim
{\bf 4.8. Proposition} The canonical isomorphisms $\gr\A(n)\to\P(n)$
define isomorphisms $\gr\A_m\to\P_m$ and $\gr\A\to\P$.
\endproclaim

\Proof The proof is the same as that of Proposition 2.8. $\square$
 
\bigskip
\noindent
{\bf 4.9.} Let us identify $\P(n)$ with the algebra of 
polynomial functions in the matrix elements of 
a matrix $x=(x_{ij})_{i,j=-n}^n$
such that $x^t=-x$, where $(x^t)_{ij}=\theta_{ij}x_{-j,-i}$; see 3.1.
\proclaim
{\bf Proposition} Any element $\phi$ of the algebra $\P_m(n)$ 
such that $\deg\phi<n-m$ can be represented as a polynomial
in the functions  
$$
p_n^{(M)}(x)=\tr(x^M),\tag1
$$
and
$$
p_{ij|n}^{(M)}(x)=(x^M)_{ij}\tag2
$$
with $|i|,|j|\leq m$.
\endproclaim

\Proof Let us introduce submatrices $a,b,c,d$ of the matrix $x$
according to the decomposition $n=m+(n-m)$. Namely, $a$ and $d$
are the square matrices whose rows and columns are enumerated by 
the indices $i,j\in\{-m,\dots,m\}$ and $i,j\notin\{-m,\dots,m\}$,
respectively, while $b$ corresponds to the subsets of indices
$i\in\{-m,\dots,m\}$ and $j\notin\{-m,\dots,m\}$, and $c=-b^t$.
(In other words, to obtain $a,b,c,d$ we represent the set of indices
$\{-n,\dots,n\}$ as the union
$\{-m,\dots,m\}\cup\{-n,\dots,-m-1,m+1,\dots,n\}$
and take the corresponding decomposition of the matrix $x$ in blocks).

We have to prove that any
polynomial $\phi(x)$ of degree $<n-m$ 
satisfying the invariance condition
$$
\phi(x)\equiv\phi(a,b,c,d)=\phi(a,bg^{-1},gc,gdg^{-1}),
\quad g\in \G_m(n)
$$
is a polynomial function of the invariants $\tr(x^M)$ and $(x^M)_{ij}$;
here we regard elements $g\in \G_m(n)\simeq \G(n-m)$ as matrices
of the same shape as $d$.
Repeating the arguments of the proof of Proposition 2.9 we
see that these invariants can be replaced by invariants of 
the form $a_{ij}$,
$\tr(d^M)$ with $M\geq 1$, and $(bd^{M-2}c)_{ij}$ with $M\geq 2$.

Since $b=-c^t$ our task is reduced to
describing the polynomial invariants
$\phi(c,d)$, where $c$ lies in the direct sum of
several copies of the $\G_m(n)$-module $V$
of column vectors of the length $2(n-m)$, and $d$ is
an element of the adjoint module $\g_m(n)$.
Any such invariant can be decomposed into a sum of expressions of the 
form
$$
\psi(c,\dots,c,d,\dots,d),
$$
where $\psi$ is a multilinear invariant. That is, $\psi$ is a morphism
of $\G_m(n)$-modules
$$
\psi:\ \undersetbrace q\to{V\otimes\dots\otimes V}
\otimes\undersetbrace r\to{\g_m(n)\otimes\dots\otimes\g_m(n)}
\to\C, \tag3
$$
for some $q$ and $r$, where we consider $\C$ as the trivial module. Note
that  the adjoint module $\g_m(n)$ is isomorphic to $\Lambda^2(V)$ in
the orthogonal case, and to $S^2(V)$ in the symplectic case.
Let us identify the exterior or symmetric square of the module $V$
with a submodule in the tensor product $V\otimes V$. Then any morphism
(3) can be obtained as the restriction of a morphism 
of $\G_m(n)$-modules of the form
$$
\chi:\ \undersetbrace q \to{V\otimes\dots\otimes V}\otimes
\undersetbrace {2r} \to{V\otimes\dots\otimes V}\to \C. 
$$
Consider the invariant nondegenerate
bilinear form on the space $V$ which is preserved by the action of
$\G_m(n)$:
$$
\alignat2
&\langle e_i,e_j\rangle=\delta_{i,-j}&&\quad\text{in the orthogonal
case},\\
&\langle e_i,e_j\rangle=\sgn i\cdot\delta_{i,-j}&&\quad
\text{in the symplectic
case},
\endalignat
$$
where $\{e_i\}$ is the canonical basis of $V\simeq\C^{2(n-m)}$,
$i=-n,\dots,-m-1,m+1,\dots,n$.
Using the classical invariant theory for the orthogonal and
symplectic groups (see H. Weyl [W], Theorems 2.9.A and 6.1.A),
we obtain that the invariant
$$
\chi(c_1,\dots,c_q,u_1,v_1,\dots,u_r,v_r),\quad c_i,u_i,v_i\in V
$$
can be represented as a sum of monomials
in bilinear invariants $\langle c_i,c_j\rangle$, $\langle c_i,u_j\rangle$,
$\langle c_i,v_j\rangle$, $\langle u_i,u_j\rangle$,
$\langle u_i,v_j\rangle$, and $\langle v_i,v_j\rangle$
such that in
each of the monomials each letter appears exactly once. 
Here we have used the assumption $\deg\phi<n-m$,
which implies that $q+r<n-m$ and hence $q+2r<2(n-m)$.
This has enabled us to eliminate the invariants of the form
$\det\ts [w_1,\dots,w_{2(n-m)}],\ w_i\in V$, which could occur
in the orthogonal case. (In the symplectic case
this assumption is not essential).
Note that if for some $i$ we permute the letters $u_i$ and $v_i$
in such a monomial, then its restriction will still determine, up to a
sign,
the same invariant of the form (3). Therefore we may only consider
the monomials which are products of submonomials of the form
$$
\langle v_{s_1},u_{s_2}\rangle\langle v_{s_2},u_{s_3}\rangle\dots
\langle v_{s_k},u_{s_1}\rangle \tag4
$$
and
$$
\langle c_{\alpha},u_{s_1}\rangle
\langle v_{s_1},u_{s_2}\rangle\dots
\langle v_{s_k},c_{\beta}\rangle\qquad
\text{(in particular,} \quad 
\langle c_{\alpha},c_{\beta}\rangle).
\tag5
$$
Let us check that the monomials (4) and (5) determine
the invariants of the form $\tr(d^k)$ and $(c^td^kc)_{ij}$, respectively.

Consider first the orthogonal case and calculate $\chi(d)$, where
the invariant $\chi$ is determined by a monomial of the form (4).
We can obviously assume that $(s_1,\dots,s_k)=(1,\dots,k)$.
That is,
$$
\chi:\ V^{\ot 2k}\to \C,\qquad
\chi:\ u_1\ot v_1\ot\cdots\ot u_k\ot v_k
\mapsto
\langle v_{1},u_{2}\rangle\langle v_{2},u_{3}\rangle\dots
\langle v_{k},u_{1}\rangle.
\tag 6
$$
We shall regard the matrix $d$ as an element of $V\ot V$ 
using the vector space isomorphism 
$\g_m(n)\to\Lambda^2(V)$ defined by the formula
$
F_{ij}\mapsto e_i\otimes e_{-j}-e_{-j}\otimes e_{i}.
$
We have $d=\sum d_{ij}E_{ij}={1\over2}\sum d_{ij}F_{ij}$, because
$d_{ij}=-d_{-j,-i}$. Therefore,
$$
\align
\undersetbrace k \to{d\otimes\dots\otimes d}&{}=
{1\over2^k}\sum d_{i_1j_1}\dots d_{i_kj_k}F_{i_1j_1}\otimes\dots\otimes
F_{i_kj_k}\mapsto \\
{1\over2^k}\sum d_{i_1j_1}\dots d_{i_kj_k}
&(e_{i_1}\ot e_{-j_1}- e_{-j_1}\ot e_{i_1})\ot\dots\ot
(e_{i_k}\ot e_{-j_k}- e_{-j_k}\ot e_{i_k})\\
{}={}&\sum d_{i_1j_1}\dots d_{i_kj_k}
(e_{i_1}\ot e_{-j_1})\ot\dots\ot(e_{i_k}\ot e_{-j_k}).
\endalign
$$
So, by (6)
$$
\align
\chi(d)={}&
\sum d_{i_1j_1}\dots d_{i_kj_k}
\langle e_{-j_1},e_{i_2}\rangle\langle e_{-j_2},e_{i_3}\rangle
\cdots\langle e_{-j_k},e_{i_1}\rangle\\
{}={}&\sum d_{i_1j_1}\dots d_{i_kj_k}
\delta_{j_{1}i_{2}}\delta_{j_{2}i_{3}}\cdots \delta_{j_{k}i_{1}}=
\sum d_{i_{1}i_{2}}d_{i_{2}i_{3}}\dots d_{i_{k}i_{1}}=\tr(d^k).
\endalign
$$

In the symplectic case
the vector space
isomorphism $\g_m(n)\to S^2(V)$ can be defined by the formula
$
F_{ij}\mapsto\sgn j\ts(e_i\otimes e_{-j}+ e_{-j}\otimes e_{i}).
$
An analogous calculation shows that here for an invariant $\chi$ 
of the form (4) one has
$
\chi(d)=(-1)^k\tr(d^k).
$

For the monomials of the form (5) the calculation is quite similar and
will be omitted. $\square$

\bigskip
\noindent
{\bf 4.10.} Since $x^t=-x$, it follows from (4.9.1) and (4.9.2) that
$$
p_n^{(M)}(x)=(-1)^M p_n^{(M)}(x),\qquad 
p_{ij|n}^{(M)}(x)=(-1)^M \theta_{ij}p_{-j,-i|n}^{(M)}(x).
$$
So, $p_n^{(M)}$ vanishes for odd $M$, whereas for the elements
$p_{ij|n}^{(M)}$ we may impose the following restrictions on $i,j,M$:
\smallskip

in the orthogonal case:
$
\ \ 
i+j<0\quad\text{for}\ \  M \ \ \text{odd,}\qquad
i+j\leq 0\quad\text{for}\ \  M \ \ \text{even};
$

in the symplectic case:
$
\ \ 
i+j<0\quad \text{for}\ \  M \ \ \text{even,}\qquad
i+j\leq 0\quad\text{for}\ \  M \ \ \text{odd.}
$
\smallskip

Fix $K\in\{1,2,\dots\}$, and 
$m\in\{0,1,\dots\}$ in the case of
$\g(n)=\oa(2n),\ \spa(2n)$ and $m\in\{-1,0,1,\dots\}$ in the case of
$\g(n)=\oa(2n+1)$. 

\medskip
\proclaim
{\bf Proposition} Let $|i|,|j|\leq m$ and $1\leq M\leq K$.
Suppose $n$ is large enough. Then the elements
$p_n^{(M)}$ with $M$ even and the elements
$p_{ij|n}^{(M)}$, where $i,j,M$ satisfy the above restrictions,
are algebraically independent.
\endproclaim

\Proof We shall only consider the orthogonal case. The proof
in the symplectic case can be obtained by an obvious adjustment.

For any triple $(i,j,M)$ 
satisfying our assumptions we choose
a subset 
$$
\Omega_{ijM}\subset\{m+1, m+2,\dots\}
$$
of cardinality $M-1$ in such a way that all these subsets be disjoint.
Let $n$ be so large that all of them belong to $\{m+1, m+2,\dots,n-K\}$.
Introduce complex parameters 
$$
y_1,\dots y_K \quad\text{and}\quad z_{ijM},
$$
where $|i|,|j|\leq m$ and $M=1,\dots K$.

Let us now define a linear
operator $x_{ijM}$ in $\C^N$ depending on $z_{ijM}$ as follows.
Let $e_{-n},\dots,e_{n}$ be the canonical basis in $\C^N$ and
$a_1<\dots<a_{M-1}$ be all the  elements of $\Omega_{ijM}$. Then
for $i+j<0$ 
$$
\align
x_{ijM}:\quad &e_j\mapsto z_{ijM}e_{a_{M-1}},\ \ 
e_{a_{M-1}}\mapsto e_{a_{M-2}},\ \dots, \ e_{a_1}\mapsto e_i,\\
&e_{-i}\mapsto -e_{-a_1},\ \  e_{-a_1}\mapsto -e_{-a_2},\
\dots, \ e_{-a_{M-1}}\mapsto -z_{ijM}e_{-j},\\
&e_k\mapsto 0\quad
\text{for}\quad  k\notin \{-i\}\cup\{j\}\cup \pm\Omega_{ijM};
\endalign
$$
while for $i+j=0$
$$
\align
x_{i,-i,M}:\quad &e_{-i}\mapsto
z_{i,-i,M}e_{a_{M-1}}-\tfrac12 e_{-a_1},\ \ 
e_{a_{M-1}}\mapsto e_{a_{M-2}},\ \dots,
\ e_{a_1}\mapsto \tfrac12 e_i,\\
&e_{-a_1}\mapsto -e_{-a_2},\
\dots, \ e_{-a_{M-1}}\mapsto -z_{i,-i,M}e_i,\\
&e_k\mapsto 0\quad
\text{for}\quad  k\notin \{\pm i\}\cup\pm\Omega_{i,-i,M}.
\endalign
$$

Note that $(x_{ijM})^t=-x_{ijM}$. Moreover, for the $(i,-i)$-matrix
element of the $M$th power of the operator $x_{i,-i,M}$ we have
$$
\left((x_{i,-i,M})^M\right)_{i,-i}=
\cases z_{i,-i,M}, \quad&\text{if}\quad M \quad\text{even},\\
0, \quad&\text{otherwise}.
\endcases
$$

Now define a linear operator $x$ in $\C^N$ depending on all 
the variables by setting
\medskip
$$
xe_a=\cases \sum_{i,j,M}  x_{ijM}e_a,\quad &\text{for}\quad |a|\leq
n-K,\\ \sgn a\cdot y_{|a|-(n-K)}e_a,\quad &\text{for}\quad n-K<|a|\leq n.
\endcases
\tag1 
$$
\medskip
\noindent
Then for any matrix $x$ of the form (1) we have
$$
\align
p_n^{(M)}(x)={}&2(y_1^M+\dots+y_K^M)+\phi(\ldots,\  z_{abL},\ldots);\\
p_{ij|n}^{(M)}(x)={}&z_{ijM}+\psi(\ldots,\  z_{abL},\ldots),\qquad
L<M; 
\endalign
$$
where the indices $i,j,M$ satisfy the assumptions of the proposition
($\phi$ and $\psi$ do not depend on $y_1,\dots,y_K$).
Thus our polynomials are algebraically independent even if
they are restricted to the affine subspace
of matrices of the form (1). $\square$

\bigskip
\noindent
{\bf 4.11.} If follows from (4.9.1) and (4.9.2)
that
$$
p_n^{(M)}\in p_{n-1}^{(M)}+\I'(n),\qquad
p_{ij|n}^{(M)}\in p_{ij|n-1}^{(M)}+\I'(n);
$$
cf. 2.11. Hence, the sequence
$$
p^{(M)}:=(p_n^{(M)}|\ n\geq 0)
$$
is a well-defined element of the algebra $\P_0$ in the case of
$\g(n)=\oa(2n),\ \spa(2n)$ and of the algebra $\P_{-1}$ in the case of
$\g(n)=\oa(2n+1)$,
while the sequence
$$
p_{ij}^{(M)}:=(p_{ij|n}^{(M)}|\ n\geq m)
$$
is a well-defined element of the algebra $\P_m$.

\proclaim
{\bf Proposition} The following elements are algebraically
independent generators of the commutative algebra $\P_m$
{\rm(}below $|i|,|j|\leq m${\rm)}.

{\rm (i)} In the orthogonal case: 
$$
\align
&p^{(M)},\quad\text{where}\quad M=2,4,6,\dots,\\
&p_{ij}^{(M)},\quad\text{where}\quad i+j<0\quad\text{for}\ \  M
\ \ \text{odd},\quad\text{and}\quad i+j\leq 0\quad\text{for}\ \  M
\ \ \text{even.}
\endalign
$$

{\rm (ii)} In the symplectic case: 
$$
\align
&p^{(M)},\quad\text{where}\quad M=2,4,6,\dots,\\
&p_{ij}^{(M)},\quad\text{where}\quad i+j<0\quad\text{for}\ \  M
\ \ \text{even},\quad\text{and}\quad i+j\leq 0\quad\text{for}\ \  M
\ \ \text{odd.}
\endalign
$$

Dropping the restriction $|i|,|j|\leq m$, we obtain algebraically
independent generators
of the commutative algebra $\P$.
\endproclaim

\Proof Let $p=(p_n|\ n\geq m)\in\P_m$. By Proposition
4.9, for any $n$ such that $\deg p<n-m$
the element $p_n\in \P_m(n)$
can be represented as a polynomial $\phi_n$ in the variables $p_n^{(M)},\
p_{ij|n}^{(M)}$ with $|i|,|j|\leq m,\ M\leq \deg p$. By Proposition 4.10,
$\phi_n$ does not  depend on $n$ for sufficiently large $n$. Hence, the
elements  $p^{(M)}$ and $p_{ij}^{(M)}$ are generators. Their algebraic
independence is clear from Proposition 4.10. $\square$

\medskip
Next few subsections contain preliminaries for the proof of the
main theorem of this section, Theorem 4.17.

\bigskip
\noindent
{\bf 4.12.} Let us recall the definition of the {\it twisted Yangian\/}
(see [MNO], Section 3 for further details). 
This is an associative algebra naturally 
associated with
the orthogonal or symplectic Lie algebra and denoted by
$\Y^+(2n)$, $\Y^+(2n+1)$, and $\Y^-(2n)$ for the cases of
$\oa(2n)$, $\oa(2n+1)$, and $\spa(2n)$, respectively. We shall
consider all the three cases simultaneously. Consider the Yangian
$\Y(N)$ (see 2.14), where we assume that the generators $t_{ij}^{(M)}$
are enumerated by the indices 
$i,j$ running
through the set $\{-n,\dots,-1,0,1,\dots n\}$ if $N=2n+1$ and
through the set $\{-n,\dots,-1,1,\dots n\}$ if $N=2n$.

Let us introduce the $S$-{\it matrix\/} $S(u)=(s_{ij}(u))$ by setting
$
S(u):=T(u)T^t(-u),
$
or, in terms of matrix elements,
$$
s_{ij}(u)=\sum_a \theta_{aj}t_{ia}(u)t_{-j,-a}(-u). 
$$
Write
$$
s_{ij}(u)=\delta_{ij}+s_{ij}^{(1)}u^{-1}+s_{ij}^{(2)}u^{-2}+\cdots.
$$
The twisted Yangian $\Y^{\pm}(N)$ is the subalgebra of $\Y(N)$ generated by
the elements $s_{ij}^{(1)},s_{ij}^{(2)},\dots$, where $-n\leq i,j\leq n$.

One can show that the $S$-matrix satisfies the
`quaternary relation' (0.3) and the `symmetry relation' (0.4);
see [MNO], Section 3 for the proof. 
They can be rewritten in terms of 
the generating series $s_{ij}(u)$ as follows:
$$
\aligned
[s_{ij}(u),s_{kl}(v)]={}&{1\over
u-v}(s_{kj}(u)s_{il}(v)-s_{kj}(v)s_{il}(u))\\
-{}&{1\over u+v}(\theta_{k,-j}s_{i,-k}(u)s_{-j,l}(v)-
\theta_{i,-l}s_{k,-i}(v)s_{-l,j}(u))\\
+{}&{1\over u^2-v^2}(\theta_{i,-j}s_{k,-i}(u)s_{-j,l}(v)-
\theta_{i,-j}s_{k,-i}(v)s_{-j,l}(u)) 
\endaligned
\tag 1
$$
and
$$
\theta_{ij}s_{-j,-i}(-u)=s_{ij}(u)\pm {s_{ij}(u)-s_{ij}(-u)\over 2u}. 
\tag 2
$$
These are defining relations for the twisted Yangian $\Y^{\pm}(N)$, so that
one can regard $\Y^{\pm}(N)$ either as a subalgebra in $\Y(N)$ or
as an abstract algebra with the generators $s_{ij}^{(M)}$ and the 
relations (1) and (2).

There is an analog of the quantum determinant for the
twisted Yangians. It is denoted by $\sdet S(u)$ and is called the
{\it Sklyanin determinant\/}. This is a formal series in $u^{-1}$
and its coefficients generate the center of the algebra
$\Y^{\pm}(N)$. The Sklyanin determinant is related
to the quantum determinant by the formula:
$$
\sdet S(u)=\gamma_n(u)\ts\qdet T(u)\ts\qdet T(-u+N-1),
$$
where 
$$
\gamma_n(u)=\cases 1\qquad&\text{for}\quad \Y^+(N),\\
\dfrac {2u+1}{2u-2n+1} \qquad&\text{for}\quad \Y^-(2n).
\endcases
\tag 3
$$
An explicit determinant-type expression analogous to (2.14.5)
for $\sdet S(u)$ in terms of the generators $s_{ij}(u)$ was given in [M2].

\bigskip
\noindent
{\bf 4.13.} Denote by $F$ the $N\times N$-matrix
formed by the elements $F_{ij},\ -n\leq i,j\leq n$,
$$
F=\sum_{i,j}F_{ij}\otimes E_{ij}\in\A(n)\otimes\End(\C^N).
$$
It was proved in [MNO], Proposition 3.11 that the mapping
$$
\eta:\ S(u)\mapsto 1+\frac{F}{u\pm 1/2}
\tag 1
$$
defines an algebra homomorphism
$
\Y^{\pm}(N)\to \A(n). 
$
This implies that the coefficients of the series $\eta(\sdet S(u))$
belong to the center $\Z(n)$ of the algebra $\A(n)$.
The images of these coefficients under the Harish-Chandra isomorphism
$\omega:\ \Z(n)\to\M(n)$
(see 1.8) were found in different ways in [M2], Section~5 and 
[MN], Section~6. The result can be written as follows:
$$
\omega:\ \gamma_n(u)\ts\eta(\sdet S(-u+\frac N2-1))\mapsto
\prod_{i=1}^n \frac{(u+1/2)^2-l_i^2}{(u+1/2)^2-\rho_i^2},
\tag 2
$$
where $\gamma_n(u)$ is defined in (4.12.3).

\bigskip
\noindent
{\bf 4.14.} Introduce the following series in $u^{-1}$ with coefficients
in $\M(n)$:
$$
\chi_n(u)=\frac{u+\rho_1-c+1/2}{u+\rho_1+1/2}\cdot
\prod_{i=1}^n \frac{(u+1/2)^2-l_i^2}{(u+1/2)^2-(\rho_i-c)^2}. \tag 1
$$
Note that $(\rho_i-c)^2$ in the denominator is equal to
$(l_i^c)^2$. Therefore, 
using definition (3.7.1) we immediately obtain
that 
$$
\pi_{n,c}: \chi_n(u)\mapsto \chi_{n-1}(u), \tag 2
$$
and so, the sequence $\chi(u)=(\chi_n(u))$ 
is an element of $\Z_c[[u^{-1}]]$. This fact will be used later.

In the next proposition we regard
$\chi_n(u)$ as an element of $\Z(n)[[u^{-1}]]$
identifying $\M(n)$ with $\Z(n)$ via the Harish-Chandra isomorphism
$\omega$. Set 
$$
\kappa_n:=\frac{N\mp1}{2}. \tag 3
$$

\medskip
\proclaim
{\bf Proposition} The mapping
$$
\varphi_n: S(u)\mapsto {u+c+\kappa_n\over u+\kappa_n}\ts\chi_n(u)
\left(1-{F\over u+\kappa_n}\right)^{-1}
\tag 4
$$
defines an algebra homomorphism
$$
\varphi_n: \Y^{\pm}(N)\to\A(n). \tag 5
$$
\endproclaim

\Proof Define the matrix
$\widehat S(u)$ by the
following formula:
$$
\sdet S(u)=\widehat S(u)S(u-N+1).\tag 6
$$
It was proved in [M4], Proposition 1.1 that the mapping
$$
S(u)\mapsto \gamma_n(u)\ts \wh S(-u+\frac N2-1)\tag 7
$$
defines an automorphism of the algebra $\Y^{\pm}(N)$.
Furthermore, it is obvious from the defining relations (4.12.1) and
(4.12.2) 
that any
invertible even formal series $g(u)\in\C[[u^{-2}]]$
also defines an automorphism of the algebra $\Y^{\pm}(N)$
given by
$
S(u)\mapsto g(u)S(u).
$
Fix the following series $g(u)=g_n(u)$,
$$
g_n(u)=
\prod_{i=1}^{n-1}\frac{u^2-(\rho_{i}-1/2)^2}{u^2-(\rho_{i}-c-1/2)^2}.
$$
Taking the composition of these automorphisms with the homomorphism $\eta$
(see (4.13.1)) we obtain another algebra homomorphism
$
\Y^{\pm}(N)\to \A(n) 
$
such that 
$$
S(u)\mapsto g_n(u)
\gamma_n(u)\ts \eta(\wh S(-u+\frac N2-1)).
$$
Let us check that this homomorphism coincides with $\varphi_n$.
Indeed, by definition (6),
$$
\wh S(-u+\frac N2-1)=\sdet S(-u+\frac N2-1)\ts
\left(S(-u-\frac N2)\right)^{-1}.
$$
We have
$$
\eta:\ \left(S(-u-\frac N2)\right)^{-1}\mapsto 
\left(1-{F\over u+\kappa_n}\right)^{-1}.
$$
A direct calculation with the use of (4.13.2) 
and the relations $\rho_{i+1}=\rho_i-1$, $\kappa_n=-\rho_n+1/2$
shows that
$$
g_n(u)\ts \gamma_n(u)\ts \eta(\sdet S(-u+\frac N2-1)) =
{u+c+\kappa_n\over u+\kappa_n}\ts\chi_n(u). \qquad\square
$$
 
\bigskip
\noindent
{\bf 4.15.} If follows from the defining relations (4.12.1), (4.12.2)
and the
Poincar\'e--Birkhoff--Witt theorem for the twisted Yangian 
(see [MNO], Remark 3.14) that for any $N\geq 3$
one has natural inclusions
$$
\Y^{\pm}(N-2)\hra\Y^{\pm}(N). \tag1
$$
Assume that $0\leq m\leq n$
in the case of
$\g=\oa(2n),\spa(2n)$, and $-1\leq m\leq n$
in the case of
$\g=\oa(2n+1)$ and set $M=2m$ and $M=2m+1$,
respectively. 
Then using (1) we can regard $\Y^{\pm}(M)$
as a subalgebra in $\Y^{\pm}(N)$.

\proclaim
{\bf Proposition} The image of the restriction of the homomorphism
$\varphi_n$ to the subalgebra $\Y^{\pm}(M)$ is contained
in the centralizer $\A_m(n)$.
\endproclaim

\Proof It follows from (4.12.1)
that
$$
[s_{kl}^{(1)}, s_{ij}(u)]=\delta_{il}s_{kj}(u)-\delta_{kj}s_{il}(u)
-\theta_{i,-l}\delta_{k,-i}s_{-l,j}(u)+
\theta_{k,-j}\delta_{-j,l}s_{i,-k}(u).\tag2
$$
In particular,
$$
[s_{kl}^{(1)}, s_{ij}(u)]=0 \quad\text{for}\quad
|i|,|j|\leq m<|k|,|l|. \tag3
$$
On the other hand, it is easy to verify that the image of $s_{kl}^{(1)}$
under the homomorphism $\varphi_n$ coincides with
$F_{kl}$. So, (3) implies that
$
[F_{kl}, \varphi_n(s_{ij}(u))]=0.
$
$\square$

\bigskip
\proclaim
{\bf 4.16. Proposition} The sequence of homomorphisms 
$(\varphi_n|\ n\geq m)$ defines a homomorphism
$$
\varphi: \Y^{\pm}(M)\to\A_m. \tag1
$$
\endproclaim

\Proof We have to verify that the homomorphisms $\varphi_n$ 
are compatible with the sequence of morphisms (4.4.1), that is, the
following diagram is commutative:
$$
\CD
\Y^{\pm}(M) @= \Y^{\pm}(M) @= \cdots @= \Y^{\pm}(M)@=\cdots\\
@V \varphi_m VV     @V\varphi_{m+1}VV @. @V\varphi_{n}VV\\
\A_m(m) @<<\pi_{m+1,c}< \A_{m}(m+1)@<<<\cdots 
@<<\pi_{n,c}< \A_m(n)@<<<\cdots.
\endCD
$$

Let us set
$$
\Sym(u)={u+c+\kappa_n\over u+\kappa_n}
\left(1-{F\over u+\kappa_n}\right)^{-1} \tag 2
$$
and denote the matrix elements of $\Sym(u)$ by $\sigma_{ij|n}(u)$.
We need to prove that
$$
\chi_n(u)\sigma_{ij|n}(u)-\chi_{n-1}(u)\sigma_{ij|n-1}(u)
\in \I(n)[[u^{-1}]],\qquad |i|,|j|\leq n-1. \tag 3
$$
However, by (4.14.2)
$$
\chi_n(u)-\chi_{n-1}(u)
\in \I(n)[[u^{-1}]]. \tag 4
$$
Hence, since the coefficients of $\chi_n(u)$ belong to $\M(n)\simeq\Z(n)$,
to prove (3) we only need to show that
$$
\sigma_{ij|n}(u)-\sigma_{ij|n-1}(u)
\in \I(n)[[u^{-1}]],\qquad |i|,|j|\leq n-1. \tag 5
$$
We shall prove by induction on $k$ that for the coefficients
of the series $\sigma_{ij|n}(u)$ one has:
$$
\alignedat2
\text{(i)}\quad &\sigma_{in|n}^{(k)}\in\I(n),
\quad &&-n\leq i\leq n,\quad k\geq 1;\\
\text{(ii)}\quad &\sigma_{ij|n}^{(k)}-\sigma_{ij|n-1}^{(k)}\in\I(n),
\quad &&-n+1\leq i,j\leq n-1,\quad k\geq 1.
\endalignedat\tag6
$$
By definition of $\Sym(u)$,
$$
\Sym(u)(u+\kappa_n-F)=u+c+\kappa_n.
$$
Hence
$$
u\ts\Sym(u)=u+c+\kappa_n+\Sym(u)(F-\kappa_n).\tag7
$$
Write $\Sym(u)=\Sym^{(0)}+\Sym^{(1)}u^{-1}+\dots$\ts.
Then (7) implies that $\Sym^{(0)}=1$ and
$$
\align
\Sym^{(1)}&=F+c,\tag8\\
\Sym^{(k)}&=\Sym^{(k-1)}(F-\kappa_n),
\quad k\geq 2.\tag9
\endalign
$$
By (8) we have
$$
\sigma_{in|n}^{(1)}
=F_{in}+\delta_{in}c\in\I(n),\quad -n\leq i\leq
n, 
$$
and
$$
\sigma_{ij|n}^{(1)}-\sigma_{ij|n-1}^{(1)}=0,
\quad -n+1\leq i,j\leq n-1.
$$ 
So, we have verified (6) for $k=1$. For $k>1$ we obtain from (9)
that
$$
\align
\sigma_{in|n}^{(k)}=
&\sum_{a=-n}^n\sigma_{ia|n}^{(k-1)}(F_{an}-
\delta_{an}\kappa_n)\\
=&\sum_{a=-n}^{n-1}\sigma_{ia|n}^{(k-1)}F_{an}+
\sigma_{in|n}^{(k-1)}(F_{nn}+c)-(c+\kappa_n)
\sigma_{in|n}^{(k-1)}.
\endalign
$$
By the induction hypotheses, this expression lies in $\I(n)$, which
proves (i) in (6). Again using (9) we obtain for 
$-n+1\leq i,j\leq n-1$ that
$$
\sigma_{ij|n}^{(k)}-\sigma_{ij|n-1}^{(k)}=
\sum_{a=-n}^{n}\sigma_{ia|n}^{(k-1)}(F_{aj}-\delta_{aj}\kappa_n)
-\sum_{a=-n+1}^{n-1}\sigma_{ia|n-1}^{(k-1)}(F_{aj}-
\delta_{aj}\kappa_{n-1}),
$$
which can be rewritten as
$$
\align
\sigma_{ij|n}^{(k)}-&\sigma_{ij|n-1}^{(k)}=\\
&\sum_{a=-n+1}^{n-1}(\sigma_{ia|n}^{(k-1)}-
\sigma_{ia|n-1}^{(k-1)})F_{aj}\tag10\\
&-\kappa_n\sigma_{ij|n}^{(k-1)}+\kappa_{n-1}\sigma_{ij|n-1}^{(k-1)}
+\sigma_{in|n}^{(k-1)}F_{nj}\tag11\\
&+\sigma_{i,-n|n}^{(k-1)}F_{-n,j}.\tag12
\endalign
$$
Let us consider the summands on the right hand side separately.
By Proposition 4.14 the relations (4.12.1) are satisfied by
the matrix elements of the matrix $\chi_n(u)\Sym(u)$. Then this is also
true
for the matrix elements of $\Sym(u)$ because the coefficients
of $\chi_n(u)$ are central in $\A(n)$.
So, using
(4.15.2) and (8) we obtain for a summand in (10):
$$
\align
(\sigma_{ia|n}^{(k-1)}-{}&\sigma_{ia|n-1}^{(k-1)})F_{aj}
=F_{aj}(\sigma_{ia|n}^{(k-1)}-\sigma_{ia|n-1}^{(k-1)})\\
{}+{}&\sigma_{ij|n}^{(k-1)}-\sigma_{ij|n-1}^{(k-1)}
-\delta_{ij}(\sigma_{aa|n}^{(k-1)}-\sigma_{aa|n-1}^{(k-1)})\\
{}-{}&\theta_{a,-a}\delta_{i,-a}(\sigma_{-a,j|n}^{(k-1)}-
\sigma_{-a,j|n-1}^{(k-1)})
+\theta_{i,-j}\delta_{-j,a}(\sigma_{a,-i|n}^{(k-1)}-
\sigma_{a,-i|n-1}^{(k-1)}).
\endalign
$$
Note that $\kappa_n-\kappa_{n-1}=1$. Hence, for (11) we have:
$$
-\kappa_n\sigma_{ij|n}^{(k-1)}+\kappa_{n-1}\sigma_{ij|n-1}^{(k-1)}
+\sigma_{in|n}^{(k-1)}F_{nj}
$$
$$
=-\kappa_{n-1}(\sigma_{ij|n}^{(k-1)}-\sigma_{ij|n-1}^{(k-1)})
+\sigma_{in|n}^{(k-1)}F_{nj}-\sigma_{ij|n}^{(k-1)}
$$
$$
=-\kappa_{n-1}(\sigma_{ij|n}^{(k-1)}-\sigma_{ij|n-1}^{(k-1)})
+F_{nj}\sigma_{in|n}^{(k-1)}-\delta_{ij}\sigma_{nn|n}^{(k-1)}.
$$
Finally, (12) looks as
$$
\sigma_{i,-n|n}^{(k-1)}F_{-n,j}=-\theta_{-j,n}
\sigma_{i,-n|n}^{(k-1)}F_{-j,n}.
$$
Using (i), we find by the induction hypotheses
that all the elements (10), (11) and (12) belong to the
ideal $\I(n)$, which
completes the proof of (6). $\square$

\bigskip
\noindent
{\bf 4.17.} The following theorem is our main result. We consider the
homomorphism $\varphi=(\varphi_n|\ n\geq m)$ defined in Proposition 4.16.
We use notation (4.14.1), (4.14.3), and (4.14.4).

\medskip
\proclaim
{\bf Theorem} The homomorphism $\varphi$
is an embedding
of the twisted Yangian $\Y^{\pm}(M)$ into the algebra $\A_m$.
Moreover, one has the decomposition
$$
\A_m=\Z_c\otimes\Y^{\pm}(M), \tag1
$$
where the twisted Yangian is identified with its image under the
embedding $\varphi$, and $\Z_c=\Z_c(\g)$ is the algebra of virtual
Laplace operators.
\endproclaim

\Proof To prove the first claim of the theorem we have to verify that
the kernel of the homomorphism $\varphi$ is trivial. Recall that
$\A_m$ is a filtered algebra with $\gr\A_m=\P_m$; see Proposition 4.8.
Let us find the highest order term
$\overline{\varphi}(s_{ij}^{(k)})$
of the sequence
$\varphi(s_{ij}^{(k)})=(\varphi_n(s_{ij}^{(k)})|\ n\geq m)$. By definition,
$$
\varphi_n(s_{ij}(u))={u+c+\kappa_n\over u+\kappa_n}\ts\chi_n(u)
\left(1-{F\over u+\kappa_n}\right)_{ij}^{-1}, \tag2
$$
where $\chi_n(u)$ is defined by formula (4.14.1).
Denote by $\overline{\varphi}_n(s_{ij}^{(k)})$ and 
$\overline{\chi}_{n}^{\ts(k)}$ the images of 
$\varphi_n(s_{ij}^{(k)})$ and $\chi_{n}^{(k)}$ in the $k$th
component of $\gr\A(n)=\P(n)$. 
Then by (2),
$$
\overline{\varphi}_n(s_{ij}^{(k)})=(F^k)_{ij}+
(F^{k-1})_{ij}\overline{\chi}_{n}^{\ts(1)}+
\dots+\overline{\chi}_{n}^{\ts(k)}
=p_{ij|n}^{(k)}+
p_{ij|n}^{(k-1)}
\overline{\chi}_{n}^{\ts(1)}+\dots+\overline{\chi}_{n}^{\ts(k)}.
$$
Hence,
$$
\overline{\varphi}(s_{ij}^{(k)})
=p_{ij}^{(k)}+
p_{ij}^{(k-1)}\overline{\chi}^{\ts(1)}+\dots+\overline{\chi}^{\ts(k)},
\tag3
$$
where $\overline{\chi}^{\ts(k)}$ 
denotes the highest order term
of the sequence $\chi^{(k)}=(\chi_n^{(k)})$.
The sequence $\chi(u)=(\chi_n(u))$ is an element of $\Z_c[[u^{-1}]]$;
see (4.14.2). So, $\overline{\chi}^{\ts(k)}$
belongs to
the algebra $\P_0$ or $\P_{-1}$ depending on whether 
$\g(n)=\oa(2n),\ts
\spa(2n)$ or $\g(n)=\oa(2n+1)$. Proposition 4.11
implies that the elements of the algebra $\P_m$
$$
p_{ij}^{(2k)},\quad i+j\leq 0; \qquad p_{ij}^{(2k-1)},\quad i+j<0;
\qquad k=1,2,\dots,
$$
in the orthogonal case, and
$$
p_{ij}^{(2k)},\quad i+j<0; \qquad p_{ij}^{(2k-1)},\quad i+j\leq 0;
\qquad k=1,2,\dots,
$$
in the symplectic case, are algebraically independent over
the subalgebra $\P_0$ or $\P_{-1}$. So, by (3), the same is true
for the elements
$$
\overline{\varphi}(s_{ij}^{(2k)}),\quad i+j\leq 0; 
\qquad \overline{\varphi}(s_{ij}^{(2k-1)}),\quad i+j<0;
\qquad k=1,2,\dots,
$$
in the orthogonal case, and
$$
\overline{\varphi}(s_{ij}^{(2k)}),\quad i+j<0; 
\qquad \overline{\varphi}(s_{ij}^{(2k-1)}),\quad i+j\leq 0;
\qquad k=1,2,\dots,
$$
in the symplectic case. Now, the
Poincar\'e--Birkhoff--Witt theorem for $\Y^{\pm}(M)$ (see [MNO],
Remark 3.14) implies that the kernel of $\varphi$ is trivial.

To prove the decomposition (1) it suffices to note
that the graded algebra
$\gr\Y^{\pm}(M)$ can be identified with the subalgebra
of $\gr\A_m=\P_m$ generated by the elements
$\overline{\varphi}(s_{ij}^{(k)}),\ k=1,2,\dots$\ts. $\square$

\bigskip
\noindent
{\bf 4.18.} Using the inclusions (4.15.1) we can define the algebra
$\Y^{\pm}(\infty)$ as the corresponding inductive limit as
$N\to\infty$. Theorem 4.17 immediately implies the following corollary.

\proclaim
{\bf Corollary} One has the isomorphism
$
\A=\Z_c\ot\Y^{\pm}(\infty). \qquad \square
$
\endproclaim
\medskip

The following are analogs of Theorem 2.21 and Proposition 2.22
for the Lie algebra $\g=
\oa(2\infty),
\spa(2\infty),\oa(2\infty+1)$. They are proved by the same argument.

\bigskip
\proclaim
{\bf 4.19. Theorem} Any $\g$-module $V\in\Omega(c)$ has a natural
$\A$-module structure.\qquad $\square$
\endproclaim

\medskip
\proclaim
{\bf 4.20. Proposition} Let $V$ be a $\g$-module with the highest weight
$\l$. Then every element $a\in\A_0$ in the case of $\g=\oa(2\infty),\ 
\spa(2\infty)$, and $a\in\A_{-1}$ in the case of $\g=\oa(2\infty+1)$
acts on $V$ as the scalar operator $f_a(\l)\cdot 1$, where $f_a\in\Mc$
corresponds to $a$ under the identification of $\A_0$ or $\A_{-1}$
with $\Mc$ from {\rm (3.8)}. \qquad $\square$
\endproclaim

\bigskip
\bigskip
\noindent
{\bf 5. Comments}
\bigskip
\noindent
{\bf 5.1.} At least, two aspects of the centralizer construction seem to be
nontrivial and rather surprising:

$\bullet$ first, the existence of the projections $\A_m(n)\to \A_m(n-1)$
which makes it possible to arrange the centralizer algebras into a
projective chain (2.3 and 4.3 above);

$\bullet$ second, the fact that the Yangians $\Y(m)$, $\Y^{\pm}(M)$ appear 
in the description of the limit centralizer algebra $\A_m$. 
Indeed, the Yangians are certain quantized (or deformed) enveloping
algebras, while in the definition of $\A_m$ there is no indication on
deformations!
\smallskip
An explanation of this phenomenon is that the commutation
relations for the classical Lie algebras can be expressed in an
$R$-matrix form (which, for the $B,C,D$ cases, involves a reflection
equation), see [MNO]. 
\medskip
\noindent
{\bf 5.2.} Thus, the centralizer construction shows 
that the Yangians $\Y(m)$ and
$\Y^{\pm}(M)$ are deeply connected with the classical Lie
algebras. One could even say that these Yangians are implicitly contained
in the enveloping algebras of the infinite rank Lie algebras of
type $A,B,C,D$. Note that unlike the twisted Yangians $\Y^{\pm}(M)$,
Drinfeld's Yangians of type $B,C,D$ are
not related to the corresponding Lie algebras in this way:
they cannot be projected onto the universal enveloping algebras
(see Drinfeld [D1]).
\medskip
\noindent
{\bf 5.3.} We believe that the centralizer construction can be applied to
certain Lie superalgebras. (About the Yangians corresponding to
the `strange' Lie superalgebras see 
Nazarov [N1], [N2]). There is a version of the centralizer 
construction for the symmetric groups; see [O1].
\medskip
\noindent
{\bf 5.4.} As was mentioned in the Introduction, irreducible
finite-dimensional representations of the centralizer algebra $\A_m(n)$
can be lifted to the corresponding Yangian $\Y(m)$ or $\Y^{\pm}(M)$. It
would be interesting to further study this construction. For instance,
to understand, especially for the $B,C,D$ case, what representations
appear as result of such a lifting.
\medskip
\noindent
{\bf 5.5.} In the present paper we have substantially exploited certain
stability effects. For instance, our main result can be viewed as
a description of the ``stable structure'' (or the ``stable part'' of
the defining relations) for the centralizer algebras. 

A somewhat different kind of stabilization is employed in the
definition of the category $\Omega(c)$ of modules over the infinite
rank classical Lie algebra $\frak g(\infty)=\injlim \frak g(n)$ of type
$B,C,D$. Note that any irreducible finite-dimensional $\frak g(n)$-module
can be viewed as a fragment of an irreducible $\frak g(\infty)$-module
which
belongs to $\Omega(c)$ with an appropriate choice of the parameter
$c$. (Note also that the same fact holds for the doubly infinite
version of $\frak{gl}$, see 2.23). This stability effect can be related to
R.~Howe's theory of reductive dual pairs [H]; especially, its
fermionic part. It would be interesting to develop this
idea and relate it to other stability effects in representation
theory; see, e.g., Benkart--Britten--Lemire~[BBL], 
Brylinski~[B], Stanley~[S1].
\medskip
\noindent
{\bf 5.6.} In [O2], Remark 2.1.20 one can find a description of a set of
generators for the virtual center $\Z_c$ (case of $\frak{gl}(\infty)$,
$c=0$; see 1.10--1.11 above) together with a description of their images
under
the isomorphism $\Z_c=\Lambda^*_c$. A somewhat different construction
is proposed in Gould--Stoilova [GS].

For the $B,C,D$ case, using an appropriate modification of the ideas
of [OO2], one can construct a linear basis in the algebra
$\Mc$, introduced in 3.8.
\medskip
The second author was supported by the Russian Foundation
for Basic Research under Grant 95-01-00814.

\bigskip
\bigskip
\noindent
{\bf References}
\bigskip

\itemitem{[BBL]}
{G. M. Benkart, D. J. Britten, and F. W. Lemire},
{\it Stability in modules for classical Lie algebras --
a constructive approach},
{Memoirs AMS}, {\bf 85} (1990), no. 430.

\itemitem{[B]}
{R. K. Brylinski}, 
{\it Stable calculus of the mixed tensor character. I},
in `S\'eminaire d'Alg\`ebre Paul Dubriel et Marie-Paule Malliavin',
Lecture Notes in Math., Vol. 1404, pp. 35--94.
Springer-Verlag, New York/Berlin, 1989. 

\itemitem{[CP1]}{V. Chari {\rm and} A. Pressley},
{\it Yangians and $R$-matrices}, {L'Enseign. Math.} {\bf
36} (1990), 267--302.

\itemitem{[CP2]}{V. Chari {\rm and} A.
Pressley}, 
{\it Fundamental representations of Yangians and singularities of
$R$-matrices}, {J. Reine Angew. Math.} {\bf
417} (1991), 87--128.

\itemitem{[C1]}{I. V. Cherednik},
{\it Factorized particles on the half-line and root systems},
{Theor. Math. Phys.}
{\bf 61} (1984), no. 1, 35--44.

\itemitem{[C2]}
{I. V. Cherednik},
{\it A new interpretation of Gelfand--Tzetlin bases}, {Duke Math. J.}
{\bf 54}
(1987),
563--577.

\itemitem{[D]}
{J. Dixmier},
{\it Alg\`ebres Enveloppantes}, 
{Gauthier-Villars, Paris},
1974.

\itemitem{[D1]}{V. G. Drinfeld},
{\it Hopf algebras and the
quantum Yang--Baxter equation}, 
{Soviet Math. Dokl.} {\bf
32} (1985), 254--258.

\itemitem{[D2]}{V. G. Drinfeld},
{\it A new realization of
Yangians and quantized affine algebras}, {Soviet Math. Dokl.} {\bf
36} (1988), 212--216.

\itemitem{[GS]}
{M. D. Gould, N. I. Stoilova}, {{\it Casimir invariants and
characteristic identities for} ${\gl}(\infty)$},
{J. Math. Phys.} {\bf 38} (1997), 4783--4793.

\itemitem{[H]}
{R. Howe},
{\it Remarks on classical invariant theory}, {Trans. AMS}
{\bf 313}
(1989),
539--570.

\itemitem{[KR]}
{A. N. Kirillov and N. Yu. Reshetikhin},
{\it Yangians, Bethe ansatz and combinatorics},
{Lett. Math. Phys.}
{\bf 12}
(1986),
199--208.

\itemitem{[KK]} {T. H. Koornwinder and V. B. Kuznetsov},
{\it Gauss hypergeometric function and quadratic $R$-matrix algebras},
{St. Petersburg Math.~J.} {\bf 6} (1994), 161--184.

\itemitem{[KS]}
{P. P. Kulish and E. K. Sklyanin},
{\it Algebraic structures related to reflection equations}, {J. Phys.}
{\bf A25}
(1992),
5963--5975.

\itemitem{[KJC]}
{V. B. Kuznetsov, M. F. J\o{rgensen}, P. L. Christiansen},
{\it New boundary conditions for integrable lattices},
{J. Phys. A} {\bf 28} (1995), 4639.

\itemitem{[M]} {I. G. Macdonald,} 
{\it Symmetric functions and Hall polynomials}, 2nd edition,
Oxford University Press, 1995.

\itemitem{[M1]}{A. Molev},
{\it Gelfand--Tsetlin bases for representations of Yangians},
Lett. Math. Phys.
{\bf 30} (1994),  53--60. 

\itemitem{[M2]}{A. Molev},
{\it 
Sklyanin determinant, Laplace operators, and characteristic identities for
classical Lie algebras}, {J. Math.
Phys.} {\bf 36} (1995), 923--945.

\itemitem{[M3]} {A. Molev}, {\it Noncommutative symmetric 
functions and Laplace operators
for classical Lie algebras}, {Lett. Math. Phys.}
{\bf 35} (1995), 135-143.

\itemitem{[M4]} {A. Molev},
{\it Finite-dimensional irreducible representations of twisted
Yangians},
Preprint CMA 047-97, Austral. Nat. University, Canberra;
q-alg/9711022.

\itemitem{[MN]} {A. Molev and M. Nazarov},
{\it Capelli identities for classical Lie algebras},
Preprint CMA 003-97, Austral. Nat. University, Canberra;
q-alg/9712021.

\itemitem{[MNO]}
{A. Molev, M. Nazarov and G. Olshanski},
{\it Yangians and classical Lie algebras}, 
Russian Math. Surveys
{\bf 51}:2
(1996),
205--282.

\itemitem{[N1]}
{M. Nazarov},
{\it Yangians of the `strange' Lie superalgebras},
in \lq Quantum Groups (P. P. Kulish, Ed.)', {Lect. Notes in Math.}
{\bf 1510},
Springer,
Berlin-Heidelberg,
1992, pp. 90--97.

\itemitem{[N2]}
{M. Nazarov},
{\it Yangian of the queer Lie superalgebra},
in `Symposium in Representation Theory', 
Yamagatta, 1992, pp. 8--15.

\itemitem{[NO]}
{M. Nazarov and G. Olshanski},
{\it Bethe subalgebras in twisted Yangians}, 
Commun. Math. Phys.
{\bf 178}
(1996),
483--506.

\itemitem{[NT1]}
{M. Nazarov and V. Tarasov},
{\it Yangians and Gelfand--Zetlin bases}, Proc. RIMS, Kyoto Univ.,
{\bf 30} (1994), 459--478.

\itemitem{[NT2]}
{M. Nazarov and V. Tarasov},
{\it Representations of Yangians with Gelfand--Zetlin bases}, 
to appear in {J. Reine Angew. Math.}; q-alg/9502008.

\itemitem{[OO1]}
A. Okounkov and G. Olshanski,
{\it Shifted Schur functions},
St.\,Petersburg Math. J.
{\bf 9}:2
(1998);
q-alg/9605042.

\itemitem{[OO2]}
A. Okounkov and G. Olshanski,
{\it Shifted Schur functions II. 
Binomial formula for characters of classical groups and applications},
in
\lq A.~A.~Kirillov Seminar on Representation Theory\rq, 
AMS,
Providence,
1997;
q-alg/9612025.

\itemitem{[O1]} {G. I. Olshanski}, 
{\it Extension of the algebra $U(g)$ for 
infinite-dimensional classical Lie algebras $g$, 
and the Yangians $Y(gl(m))$.}
{Soviet Math. Dokl.} {\bf 36}, no. 3 (1988), 569-573.

\itemitem{[O2]}{G. I. Olshanski},
{\it Representations of 
infinite-dimensional classical groups, limits of enveloping algebras, and
Yangians}, in \lq Topics in Representation Theory (A.~A.~Kirillov, Ed.)',
{Advances in Soviet Math.} {\bf 2}, AMS, Providence RI, 1991, pp.
1--66.

\itemitem{[O3]}{G. I. Olshanski},
{\it Twisted Yangians and
infinite-dimensional classical Lie algebras}, in \lq Quantum Groups 
(P.~P.~Kulish, Ed.)', {Lecture Notes in Math.} {\bf
1510}, Springer, Berlin-Heidelberg, 1992, pp. 103--120.

\itemitem{[S]}{E. K. Sklyanin},
{\it Boundary conditions for integrable quantum 
systems}, 
{J. Phys.} {\bf A21} (1988), 2375--2389.

\itemitem{[S1]}
{R. P. Stanley}, 
{{\it The stable behavior of some characters of} $SL(n, C)$},
{Linear and Multilinear Algebra} {\bf 16} (1984), 3--27.

\itemitem{[TF]}{L. A. Takhtajan {\rm and} L. D. Faddeev},
{\it Quantum inverse scattering method and the Heisenberg $XYZ$-model},
{Russian
Math. Surv.} {\bf 34} (1979), no. 5, 11--68.

\itemitem{[W]}
{H. Weyl},
{\it Classical Groups, their Invariants and Representations},
{Princeton Univ. Press},
Princeton NJ,
1946.

\itemitem{[\v{Z}]}
{D. P. \v{Z}elobenko},
{\it Compact Lie groups and their representations},
{Transl. of Math. Monographs}
{\bf 40} AMS,
Providence RI
1973.

\enddocument